\documentclass[11pt,leqno]{amsart}
\usepackage{amsmath}
\usepackage{amsfonts}
\usepackage{amssymb}
\usepackage{graphicx}
\usepackage{hyperref} 
\usepackage{a4wide} % margins
\usepackage{mathrsfs} % fancy calligraphic style
\usepackage{relsize} % big size of integrals
\usepackage{xcolor} % colors

\usepackage{graphicx}
\usepackage{hyperref}

\usepackage{enumerate}
\usepackage[shortlabels]{enumitem}

\newtheorem{theorem}{Theorem}[section]

\newtheorem{corollary}[theorem]{Corollary}

\newtheorem{example}[theorem]{Example}

\newtheorem{lemma}[theorem]{Lemma}

\newtheorem{proposition}[theorem]{Proposition}
\newtheorem{remark}[theorem]{Remark}

\newcommand{\sect}{\mathrm{Sect}}
\newcommand{\dist}{\operatorname{dist}}
\newcommand{\diam}{\operatorname{diam}}
\newcommand{\vol}{\operatorname{vol}}

\title[Isoperimetry and asymptotically nonnegative curvature]{Isoperimetric and Michael-Simon inequalities on manifolds with asymptotically nonnegative curvature}
\author[Debora Impera]{Debora Impera}
\address[Debora Impera]{Dipartimento di Scienze Matematiche ``Giuseppe Luigi Lagrange", Politecnico di Torino, Corso Duca degli Abruzzi, 24, Torino, Italy, I-10129}
\email{debora.impera@polito.it}

\author[Stefano Pigola]{Stefano Pigola}
\address[Stefano Pigola]{Dipartimento di Matematica e Applicazioni, Universit\`a di Milano Bicocca, via R. Cozzi 53, I-20126 Milano, Italy}
\email{stefano.pigola@unimib.it}

\author[Michele Rimoldi]{Michele Rimoldi}
\address[Michele Rimoldi]{Dipartimento di Scienze Matematiche ``Giuseppe Luigi Lagrange", Politecnico di Torino, Corso Duca degli Abruzzi, 24, Torino, Italy, I-10129}
\email{michele.rimoldi@polito.it}

\author[Giona Veronelli]{Giona Veronelli}
\address[Giona Veronelli]{Dipartimento di Matematica e Applicazioni, Universit\`a di Milano Bicocca, via R. Cozzi 53, I-20126 Milano, Italy}
\email{giona.veronelli@unimib.it}

\begin{document}

\begin{abstract}
We establish the validity of the isoperimetric inequality (or equivalently, an $L^1$
  Euclidean-type Sobolev inequality) on manifolds with asymptotically non-negative sectional curvature. Unlike previous results in the literature, our approach does not require the negative part of the curvature to be globally small. Furthermore, we derive a Michael--Simon inequality on manifolds whose curvature is non-negative outside a compact set.

The proofs employ the ABP method for isoperimetry, initially introduced by Cabr\'e in the Euclidean setting and subsequently extended and skillfully adapted by Brendle to the challenging context of non-negatively curved manifolds. Notably, we show that this technique can be localized to appropriate regions of the manifold. Additional key elements of the argument include the geometric structure at infinity of asymptotically non-negatively curved manifolds, their spectral properties -- which ensure the non-negativity of a Bakry-\'Emery Ricci tensor on a conformal deformation of each end -- and a result that deduces the validity of the isoperimetric inequality on the entire manifold, provided it holds outside a compact set.\end{abstract}

\maketitle

\section{Introduction}

The celebrated isoperimetric problem seeks to maximize the volume enclosed by an assigned perimeter. The problem involves two main aspects: identifying the isoperimetric regions -- that is, the sets (if they exist) that achieve the maximum -- and establishing bounds on the maximum volume that can be enclosed. This latter aspect is usually quantified through the so-called isoperimetric inequality, which states that
\begin{equation}    \label{e:isop}
    |\Omega|^{(n-1)/n} \le C |\partial \Omega|
\end{equation}
for some constant $C>0$, and for every sufficiently regular region $\Omega$. 
In the Euclidean case, or more generally on a smooth $n$-dimensional Riemannian manifold $(M^n,g)$ without boundary, it is well-known that the above isoperimetric inequality is equivalent to an $L^1$ Euclidean-type Sobolev inequality
\begin{equation}\label{e:L1Sob}
\|h\|_{n/(n-1)}\le C \|\nabla h\|_1,\qquad\forall\,h\in C^\infty_c(M),
    \end{equation}
where the $L^p$-norms and the gradient are the canonical ones induced by the Riemannian structure, \cite{FF}. In particular, either \eqref{e:isop} or \eqref{e:L1Sob} imply the validity of the $L^p$ Euclidean-type Sobolev inequalities
\begin{equation*}%\label{e:LSob}
\|h\|_{np/(n-p)}\le C_p \|\nabla h\|_p,\qquad\forall\,h\in C^\infty_c(M),\quad 1<p<n.
    \end{equation*}
Moreover, \eqref{e:L1Sob} and thus \eqref{e:isop} are also equivalent to the apparently more general version
\begin{equation}\label{e:L1Sob bdy}
            \int_\Omega |\nabla h| + \int_{\partial \Omega} h \ge C^{-1} \left(\int_\Omega h^{\frac{n}{n-1}}\right)^{\frac{n-1}{n}}
        \end{equation}
valid for any compact domain $\Omega\subset M$ with smooth boundary and any $h$ a positive smooth function on $\Omega$; see for instance \cite{brendleeichmar}. In what follows, we will refer to the validity of \eqref{e:isop}, \eqref{e:L1Sob} and \eqref{e:L1Sob bdy} interchangeably. 

It is by now well known that the isoperimetric inequalities above hold on smooth, complete, non-compact Riemannian manifolds 
$(M^n,g)$ with non-negative Ricci curvature and a positive asymptotic volume ratio
\[
\mathrm{AVR}(M):=\lim_{r\to\infty}\frac{\vol(B_r(o))}{\omega_nr^n}
\]
with $\omega_n$ the volume of the Euclidean unit ball of $\mathbb R^n$.
Classically, two main approaches have been used to prove such a result. The first relies on a celebrated formula by C. Croke \cite{croke}, which states that on any manifold
\[
\mathrm{Area}(\partial\Omega)^{n}\ge c_n \omega(\Omega)^{n+1}\vol(\Omega)^{n-1}.
\]
Here, $c_n$ is a purely dimensional explicit constant, while $\omega(\Omega)=\min_{x\in\Omega}\omega_x(\Omega)$, where $\omega_x(\Omega)$ is the ($(n-1)$-dimensional) measure of the
set of all unit vectors $v\in T_x\Omega$ such that the geodesic issuing from $x$
with tangent vector $v$ intersects $\partial\Omega$ and the geodesic segment from $x$ to its
first point of intersection with $\partial \Omega$ is length minimizing. It was observed by S.T. Yau that the non-negativity of the Ricci curvature permits to obtain a universal bound on $\omega_x(\Omega)$, hence on $\omega(\Omega)$. Indeed, assume that $\Omega$ is contained in some geodesic ball $B_r(x)$. Then, for any $k>1$, $B_{kr}(x)\setminus B_{r}(x)$ can be reached by minimizing geodesics emanating from $x$ with initial velocity in $\omega_x(\Omega)$, and the volume ratio decreases along geodesics due to volume comparison; see \cite[p.10]{Yau}.

The second approach introduced by T. Coulhon and L. Saloff-Coste permits to prove the isoperimetric inequality on manifolds with positive $\mathrm{AVR}$ which enjoy the volume doubling property and the Poincar\'e inequality; see \cite[Section 3]{CSC} for the precise definitions. As a consequence of the Bishop-Gromov comparison theorem and of the work of Buser \cite{buser}, the latter two assumptions are satisfied on manifolds with non-negative Ricci curvature. 

It is worth noting that both Yau's and Coulhon--Saloff-Coste's approaches appear to be quite sensitive to the global geometry of the manifold, making it nontrivial to localize their proofs outside a compact set of 
$M$.
 
More recently, a breakthrough improvement was proposed by S. Brendle \cite{brendle}, inspired by X. Cabr\'e's proof \cite{cabre} of the sharp isoperimetric inequality in 
 $\mathbb R^n$   via the so-called ABP method (named after A.D. Alexandrov, I.Ya. Bakelmann, and C. Pucci). In a sense, Brendle's approach follows a similar line to Yau's proof. However, instead of considering all minimizing geodesics from $\Omega$ to $B_{kr}(x)\setminus B_r(x)$, he cleverly selects a small subclass of geodesics. This selection is guided by the solution $u$ to the Neumann problem
\begin{equation}
    \label{e:Neumann intro}
\begin{cases}\Delta u = |\partial\Omega|/|\Omega|\\
\partial_\nu u =1.\end{cases}
\end{equation}
Choosing for any $x\in\Omega$ the geodesic with initial direction $\nabla u(x)$  proves sufficient to reach all points of 
$B_{r-\mathrm{diam}(\Omega)}(o)$, where $o$ is any fixed reference point in $\Omega$. Moreover, as in Yau's proof, since $\mathrm{Ric}\ge 0$, the increment of the volume measure along the geodesics grows at most as the $n$-th power of the length of the geodesic. This improvement has the major advantage of providing optimal constants and implying the sharp isoperimetric inequality on manifolds with non-negative Ricci curvature and positive  $\mathrm{AVR}$. This was previously achieved in \cite{AFM} by means of different methods, but only for the three-dimensional case.

The machinery introduced by Brendle is sufficiently flexible to address various geometric assumptions and inequalities. For instance, Brendle himself established Michael-Simon inequalities for $n$-dimensional hypersurfaces embedded in $\mathbb R^{n+m}$ \cite{brendle-Rn}, or in $(n+m)$-complete manifolds with non-negative sectional curvature \cite{brendle}. These inequalities yield the sharp Euclidean isoperimetric inequality in the case of minimal submanifolds, at least for codimensions $1$ and $2$. Subsequently, Brendle also derived a sharp logarithmic Sobolev inequality for submanifolds of the Euclidean space \cite{brendle-LogSob}, which was later extended by C. Yi and Y. Zheng to manifolds with non-negative sectional curvature \cite{YZ}. Assumptions involving intermediate Ricci curvatures have been explored in  \cite{mawu,LeeRicci}. Finally,
weighted counterparts have been obtained, where the $L^p$ norms and the volumes are computed with respect to the weighted measure $e^{-f}d\vol$ and the (non-negative) curvature assumption involves either the $k$-Bakry-\'Emery Ricci tensor \cite{Johne-arxiv} or the ($\infty$-)Bakry-\'Emery Ricci tensor assuming $f$ is bounded \cite{IRV}. 

Although less relevant to the scope of this paper, it is worth noting that a recent alternative approach to proving the isoperimetric inequality on manifolds with non-negative Ricci curvature and positive $\mathrm{AVR}$ 
has been proposed based on optimal transport methods; see, for instance, \cite{castillon,BE,BK,kris-CVPDE}. 

Introducing a small amount of negative curvature into this framework is a very delicate task. If the negative part of the curvature is small in a spectral sense, i.e.,
$\int_M G(x,y)\mathrm{Ric}_-(y)\,d\mathrm{vol}(y)\le 1/(n-2)$ for all $x\in M$ with $G$ the Green function of $M$, then the isoperimetric inequality holds in $M$ with a constant that approaches the Euclidean one as $\mathrm{Ric}_-(y)$ becomes smaller \cite{IRV}. In this scenario, one can produce a bounded positive solution to $\Delta \varphi + \lambda \mathrm{Ric}_-\varphi \le 0$, $\lambda>n-2$, on the whole $M$. As a consequence, when endowed with a suitably chosen conformal metric and weighted measure, $M$ possesses a non-negative Bakry-Émery Ricci tensor. Remarkably, the aforementioned spectral assumption is satisfied when either the sectional curvature or the Ricci curvature decays polynomially fast enough (with an additional relative volume comparison condition when the curvature assumption concerns $\mathrm{Ric}$), in a quantitative manner \cite[Theorem 1.2]{IRV}. However, even in such cases, some form of global control over the negative part $\mathrm{Ric}_-$ of the curvature across the entire manifold is still necessary.

Clearly, the constant in \eqref{e:isop} must depend on the total mass of $\mathrm{Ric}_-$ (see Remark \ref{rmk:constants}). However, it is natural to ask whether the isoperimetric inequality remains valid -- with a constant  $C$, which in general will depend on $M$ -- under suitable decay assumptions on the curvature at infinity.

In this direction, a possible approach has been first proposed by Y. Dong, H. Lin and L. Lu. Namely, they managed to prove a $L^1$-Sobolev type inequality of the form
\begin{equation}\label{e:L1Sob DLL}
            \int_\Omega |\nabla h| + \int_{\partial \Omega} h + B \int_\Omega h \ge C^{-1} \left(\int_\Omega h^{\frac{n}{n-1}}\right)^{\frac{n-1}{n}}
        \end{equation}
        on manifolds with asymptotically non-negative Ricci curvature, i.e., fixed a reference point $o\in M$,
\[
\min \mathrm{Ric} (x) \ge -\kappa(\dist(x,o)) g(x), 
\]
in the sense of quadratic forms for some positive monotone decreasing function $\kappa:[0,\infty)\to[0,\infty)$ satisfying 
\begin{equation}\label{e:b0 intro}
b_0(\kappa):=\int_0^{\infty} t \kappa(t)\,dt<\infty.    
\end{equation}
Dong-Lin-Lu's result has a number of strengths, among which are the facts that it requires no curvature assumptions beyond asymptotic curvature decay, and that the constants appearing in the inequality become sharp as the negative part of the curvature diminishes. Moreover, the result is flexible enough to be adapted to various problems (e.g., Michael-Simon and Log-Sobolev inequalities) and assumptions (e.g., weighted and intermediate Ricci curvature), as demonstrated in \cite{DLLa,DDLb,DDLc,LeeRicci}.
However, from the point of view of the isoperimetric inequality, it is somewhat lacking, as it includes the additional term 
$\|h\|_{L^1}$, which leads to a spurious isoperimetric type inequality of the form $|\partial \Omega|+B|\Omega|\ge C^{-1}|\Omega|^{1-1/n}$. This shortcoming arises from an additional term in the estimate of the distortion of the volume measure along the geodesics -- which appears in the presence of some negative curvature -- and seemingly cannot be avoided.\\

As a consequence, we adopt a different strategy based on two main ingredients. The first ingredient is an extension result of independent interest, stating that the validity of an isoperimetric inequality in the exterior of a compact set of $M$ extends—possibly with a different constant—to the entire manifold $M$.

\begin{theorem}\label{th:extension}
Let $(M^n, g)$ be a Riemannian manifold. Assume that $M$
has infinite volume and that the $L^1$-Sobolev inequality \eqref{e:L1Sob} holds for all $h\in C_c^\infty(M \setminus F)$
for some
compact $F\subset M$. Then, the same Sobolev inequality, possibly with
a different constant, holds for all $h\in C^\infty_c(M)$.
\end{theorem}

Note that the infinite volume assumption is automatically satisfied when $M$ is complete, as \eqref{e:L1Sob} implies that the volumes of the geodesic balls $B_r(o)$ of $M$ grow at least as $O(r^n)$; see \cite{Carron-Isop}, \cite[Theorem 3.1.5]{SaloffCoste-book}.  
The result analogous to Theorem \ref{th:extension} for $L^p$ Sobolev inequalities is well-known; see for instance \cite[Proposition 2.5]{carron-Duke} (for $p=2$) and \cite[Theorem 3.2]{PST} for general $p>1$.
The techniques in \cite{PST} can be adapted to $p=1$ once it is shown that the compact sets of $M$ either all have null $1$-capacity or all have strictly positive $1$-capacity. Since the standard potential theory available for $p > 1$ does not apply in this case, we adopt a more geometric approach instead (see Lemma \ref{lem:1cap}).

The second ingredient introduces one of the key novelties of our work. Specifically, we demonstrate that the ABP techniques, as developed by Brendle, can be localized to a restricted set of geodesics emanating from 
$\Omega$, provided that the volume attained by these geodesics is sufficiently large. 
As a first instance of this phenomenon, one can consider the case where the sectional curvature of 
$M$ is non-negative outside a compact set 
$K$. In this case, one can exploit the rather rigid structure of the manifolds with asymptotically non-negative sectional curvature, as notably highlighted by Abresch, Kasue, and Mashiko-Nagano-Otsuka \cite{abr,ka,mno}, and in particular the existence of a nontrivial ideal boundary $M_\infty$ which constitutes the link of the asymptotic cone of $M$. We suppose that $\Omega$ is far away from $K$, and we select among the geodesics emanating from points $x\in\Omega$ with speed $\nabla u(x)$ those which -- in a sense -- approximate the geodesic rays of $M$ which converge at infinity to a proper subset of $M_\infty$.  A sufficiently large part of these geodesics turns out to avoid the compact set $K$, so that the estimate on the volume distortion can be performed as if the curvature were non-negative. Accordingly, an isoperimetric inequality holds on $M\setminus K$ (see Remark \ref{Rmk:ToyCases}), and, similarly, a Michael-Simon inequality is satisfied on proper submanifolds of $M\setminus K$. Thanks to Theorem \ref{th:extension} these inequalities remain valid throughout the entire manifold $M$. As for the Michael-Simon inequality, more precisely, we are able to prove the following

\begin{theorem}\label{thm:nonnegsuboutcomp}
    Suppose that $(M^{n+m},g)$ has nonnegative sectional curvatures outside a compact set $K\subset M$. Let $o\in K$. Suppose that 
    $\textrm{AVR}>0$ on each end $E_j$ of $M$, i.e. there exists $\theta>0$ such that
    \[
    \limsup _{r\to\infty}\frac{\vol(B^M_r(o)\cap E_j)}{r^{n+m}}\geq \theta>0,\qquad\forall\,E_j.
    \] 
    Then, there exists $t_0>0$ such that for any compact $n$-dimensional submanifold $\Sigma\subset M\backslash B_{t_0}(o)$ with smooth boundary and any $h$ positive smooth function on $\Sigma$ it holds
    \begin{equation}\label{eq:MichaelSimon}
    \int_\Sigma \sqrt{|\nabla^\Sigma h|^2+h^2|H|^2} + \int_{\partial \Sigma} h \ge \xi \left(\int_\Sigma h^{\frac{n}{n-1}}\right)^{\frac{n-1}{n}},
    \end{equation}
    where $\xi=\xi(n,m,\theta)$ is a positive constant explicitly depending on $n,m$ and $\theta$.
    In particular, there exists a constant $C>0$ such that for any properly immersed $n$-dimensional submanifold $\Sigma$, and any positive compactly supported smooth function $h$ on $\Sigma$ it holds
    \begin{equation}\label{eq:MichaelSimonunbounded}
    \int_\Sigma \sqrt{|\nabla^\Sigma h|^2+h^2|H|^2}\ge C \left(\int_\Sigma h^{\frac{n}{n-1}}\right)^{\frac{n-1}{n}}.
    \end{equation}
\end{theorem}

In the assumptions of Theorem \ref{thm:nonnegsuboutcomp} (as well as of Theorem \ref{thm:asymp nonneg} below), there exists $R_0>0$ and a positive integer $j_0$ such that $M\setminus B_{R}(o)$ has $j_0$ connected component for all $R\ge R_0$, \cite[Theorem 3.1]{LT-harmonic}. We call an end $E_j$ of $M$, $j=1,\dots, j_0$, any of the connected components of $M\setminus B_{R_0}(o)$. Conversely to the case of nonnegative Ricci curvature previously considered by Brendle, in our setting we need to impose the volume growth assumption on each end, otherwise the result is false. A simple counterexample is given by a (two ends) manifold isometric to the disjoint union of $\mathbb R^2\setminus B_2(0)$ and $\mathbb S^1 \times [0,\infty)$ outside some compact set.
\medskip

Let us now move to the more general situation in which the sectional curvature of $(M^n,g)$ is asymptotically non-negative, i.e. \[
\min \sect (x) \ge -\kappa(\dist(x,o)), 
\]
with $\kappa$ as in \eqref{e:b0 intro}. In principle, the Ricci curvature could take negative values along all geodesics of $M$, so that when attempting to apply the ABP method, the volume distortion cannot be estimated directly without introducing an additive term, as in \eqref{e:L1Sob DLL}. To overcome this problem, we need to modify our strategy, drawing inspiration from the result proven under spectral assumptions in \cite{IRV}. Under the slightly stronger requirement that $\kappa(r)\le \kappa_0 r^{-n}$, we can prove the existence of a function $\varphi$ suitably converging to $1$ at infinity, which solves 
$\Delta \varphi + \lambda \mathrm{Ric}_-\varphi \le 0$ on $M\setminus B_{t_1}(o)$. As in \cite{CMT,IRV} we consider then on $M$ the conformal metric $\tilde g=\varphi^{2/\lambda}g$, with respect to which the
Bakry-\'Emery Ricci tensor $\mathrm{Ric} + \mathrm{Hess}f_n$ is non-negative on $M\setminus B_{t_1}(o)$, where $f_n=(n-2)\log\varphi/\lambda$. At this point, one would like to prove a weighted (equivalently unweighted, $f$ being bounded) isoperimetric inequality on $(M\setminus B_{t_1}(o),\tilde g)$ using the ABP method. As in the case in which the negative curvature is non-negative outside a compact set presented above, one takes $\Omega$ far enough from $B_{t_1}(o)$ and considers the solution $u$ to a weighted version of the Neumann problem \eqref{e:Neumann intro}. The main technical difficulty then lies in proving that there are enough geodesics of 
 $(M,\tilde g)$ emanating from points in $\Omega$ with initial speed $\nabla u$, which move toward a sufficiently large portion of $M_\infty$ without intersecting $B_{t_1}(o)$.
Indeed, an additional obstacle to overcome is that, while $(M,g)$ is asymptotically non-negatively curved, $(M,\tilde g)$ is not necessarily so.

With this strategy, we prove the last main theorem of this paper.

\begin{theorem}\label{thm:asymp nonneg}
    Let $(M^n,g)$ be 
 a complete smooth  $n$-dimensional Riemannian manifold and let $o\in M$. Suppose that
 \[
\min \mathrm{Sect} (x)\ge -\kappa_0 \dist(x,o)^{-n-\epsilon}
\]
outside some compact set $K\subset M$ for some constants $\kappa_0,\epsilon>0$. 
 Suppose that $\textrm{AVR}>0$ on any end $E_j$ of $M$, i.e. there exists $\theta>0$ such that
    \[
    \limsup _{r\to\infty}\frac{\vol(B_r(o)\cap E_j)}{r^n}\ge\theta>0,\qquad\forall\,E_j.
    \] 
    Then, there exists a constant $C=C(M)>0$ such that for any compact domain $\Omega\subset M$ with smooth boundary and any $h$ a positive smooth function on $\Omega$ it holds
    \begin{equation}\label{e:isoper asymp}
            \int_\Omega |\nabla h| + \int_{\partial \Omega} h \ge C \left(\int_\Omega h^{\frac{n}{n-1}}\right)^{\frac{n-1}{n}}.
        \end{equation}
\end{theorem}

As a special, but relevant, case of either Theorem \ref{thm:nonnegsuboutcomp} or \ref{thm:asymp nonneg}, we can state the following
\begin{corollary}\label{coro:nonnegoutcomp}
    Suppose that $(M^n,g)$ has nonnegative sectional curvatures outside a compact set $K\subset M$. Let $o\in K$. Suppose that 
    $\textrm{AVR}>0$ on each end $E_j$ of $M$, i.e. 
    \[
    \limsup _{r\to\infty}\frac{\vol(B^M_r(o)\cap E_j)}{r^{n}}>0.
    \] 
    Then, there exists a constant $C>0$ such that 
    \[|\Omega|^{(n-1)/n} \le C |\partial \Omega|
\] 
for any compact domain $\Omega \subset M$ with smooth boundary.
\end{corollary}
It is important to emphasize that, even under this stronger set of assumptions, the result is non-trivial, i.e., it can not be easily deduced from Theorem \ref{th:extension} and Brendle's work \cite{brendle}. Namely, let $M$ be a manifold whose sectional curvatures are non-negative on $M\setminus K$ for some compact set $K$. Then, in general there does not exist any non-negatively curved complete manifold $N$ without boundary such that $M\setminus K$ is isometric to a subset of $N$; see Example \ref{ex:cone}. 
\medskip

The paper is organized as follows. In Section \ref{sec:extension} we will prove the extension Theorem \ref{th:extension} and we will highlight some straightforward consequences and connections with ALE manifolds. In Section \ref{sec:preliminar} we will recall and discuss a number of preliminary results concerning manifolds with asymptotically non-negative sectional curvature, focusing in particular on their structure at infinity. Section \ref{sec:asymptotic} and Section \ref{sec:Michael-Simon} will be devoted to the proofs of Theorem \ref{thm:asymp nonneg} and \ref{thm:nonnegsuboutcomp}, respectively. Finally, in the last part of Section \ref{sec:Michael-Simon}, we also show how the techniques described above can be used to prove a logarithmic Sobolev inequality for compact submanifolds that lie outside a compact set larger than the support of the negative part of the curvature.

\subsection*{Acknowldgements.} D. Impera and M. Rimoldi are partially supported by INdAM-GNSAGA,
S. Pigola and G. Veronelli are partially supported by INdAM-GNAMPA . D. I. and M. R.
acknowledge partial support by the PRIN 2022 project “Real and Complex Manifolds: Geometry
and Holomorphic Dynamics” - 2022AP8HZ9.
\section{Extending isoperimetry inside compact sets}\label{sec:extension}

This section is devoted to the proof of Theorem \ref{th:extension} 

We first prove a characterization of manifolds admitting compact sets of null $1$-capacity; see for instance \cite{Troyanov} for the $p>1$ counterpart. It is worth noticing that the by now standard non-linear machinery used for $p>1$ is not yet available for $p=1$. Our simple geometric approach gives an alternative proof also for larger $p$.

In what follows,  $\mathrm{Cap}_1(F)$ is the infimum of $\|\nabla v\|_{L^1(M)}$ over all $v\in C^\infty_c(M)$ such that $v\ge 1$ on $F$. It is well known that, by truncation, one can reduce the problem to considering $v$ such that $0\le\chi_F \le v \le 1$. Moreover, by density, the infimum can be taken in $W^{1,1}_0(M)$ (which however may be smaller than $W^{1,1}(M)$ if the manifold is incomplete).

\begin{lemma}\label{lem:1cap}
    Let $(M,g)$ be an open  connected (possibly incomplete) smooth Riemannian manifold. Suppose that $\mathrm{Cap}_1(K)= 0$ for some compact $K\Subset M$ with non-empty interior. Then $\mathrm{Cap}_1(F)= 0$ for any other compact set $F\Subset M$ with non-empty interior.
\end{lemma}

\begin{proof}
    First, we show that for any $y\in M$ there exists a ball $B_\epsilon(y)\Subset M$ for some $0<\epsilon=\epsilon(y)$ small enough such that 
$\mathrm{Cap}_1(\overline{B_\epsilon(y)})=0$. To this end, pick a point $x\in \mathrm{int}(K)$ and a smooth embedded curve $\gamma$ from $x$ to $y$. Since the image of $\gamma$ is compact and $\gamma$ is embedded, for $\tau$ small enough the $\tau$-tubular neighborhood $B_\tau(\gamma)$ of the image of $\gamma$ is contractible and diffeomorphic through a $C^\infty$-diffeomorphism $\Psi$ to a contractible open set of $\mathbb R^n$. An explicit diffeomorphism can be constructed using Fermi coordinates along (a small prolongation) of $\gamma$. Now, choose $\epsilon,\delta>0$ small enough so that $\overline{B_\epsilon(y)}\subset B_\tau(\gamma)$ and $\overline{B_\delta(x)}\subset B_\tau(\gamma)\cap K$. There exists a $C^\infty$ diffeomorphism $\Phi:M\to M$ such that $\Phi$ is a compact perturbation of the identity and $\Phi(\overline{B_\delta(x)})=\overline{B_\epsilon(y)}$. It can be obtained by pulling back through $\Psi$ an analogous diffeomorphism in $\Psi (B_\tau (\gamma))$ which sends $\Psi(\overline{B_\delta(x)})$ onto $\Psi(\overline{B_\epsilon(y)})$ and is the identity outside a compact set. On the one hand, we have that $\mathrm{Cap_1}(\overline{B_\delta(x)})=0$ as $B_\delta(x)\subset K$. On the other hand, this latter implies that $\mathrm{Cap_1}(\overline{B_\epsilon(y)})=0$.
Indeed, $v$ is a test function for $\mathrm{Cap_1}(\overline{B_\epsilon(y)})$ if and only if $ v\circ \Psi$ is a test function for 
$\mathrm{Cap_1}(\overline{B_\delta(x)})$,
and $C_\Phi^{-1} \|\nabla v\|_{L^1(M)} \le \|\nabla (v\circ \Phi)\|_{L^1(M)}\le C_\Phi \|\nabla v\|_{L^1(M)}$, where the constant $C_\Phi$ depends only on the (finite) $C^1$-norm of $\Phi$.

To conclude the argument, take the family $\{ {B_{\epsilon(y)}(y)}\}_{y\in F}$, with $\epsilon(y)$ as above, and extract a finite subcovering of $F$, each of whose element has zero $1$-capacity. As the $1$-capacity is (finitely) subadditive, this concludes the proof.
\end{proof}

With the above lemma, we can mimic the proof given for $p>1$ in  \cite[Theorem 3.2]{PST} and \cite{Tr}. We provide the details for completeness.

\begin{proof}[Proof (of Theorem \ref{th:extension})]
    Let $\Omega$ be a precompact domain with smooth boundary such that $F\Subset \Omega$. Let also
$W_\epsilon \approx \partial\Omega \times (-\epsilon, \epsilon)$ be a bicollar neighborhood of $\partial\Omega$ such that $W_\epsilon\subset M \setminus F$, and let
$\Omega_\epsilon = \Omega \cup W_\epsilon$ and $M_\epsilon = (M \setminus \Omega) \cup W_\epsilon$. Note that, by assumption, the $L^{1}$-Sobolev inequality
with Sobolev constant $S > 0$ holds on $M_\epsilon$. Furthermore, the same $L^{1}$-Sobolev inequality,
holds for some constant $S_\epsilon > 0$ on the compact manifold with boundary $\Omega_\epsilon$; see \cite[p. 281]{He}. Now, let
$\rho\in C^\infty_c (M)$ be a cut-off function satisfying $0 \le \rho \le 1$, $\rho = 1$ on $M\setminus M_{\epsilon}$ and $\rho = 0$ on $M\setminus{\Omega _\epsilon}$.
Next, for any $v \in C^\infty_
c (M)$, write $v = \rho v+(1 - \rho) v$, and note that $\rho v \in C^\infty_c
 (\Omega_\epsilon)$ whereas
$(1 - \rho) v \in C^\infty
_c
(M_{\epsilon})$.
Therefore, we can apply the respective Sobolev inequalities and
get

\begin{align}\label{e:sobolev glued}
    \| v\|_{L^{1^\ast}(M)} &\le \|\rho v\|_{L^{1^\ast}(\Omega_\epsilon)} + \|(1-\rho)v\|_{L^{1^\ast}(M_\epsilon)}\\
    &\le S_\epsilon \|\nabla(\rho v)\|_{L^{1}(\Omega_\epsilon)} + S \|\nabla((1-\rho)v)\|_{L^{1}(M_\epsilon)}\nonumber\\ 
    &\le \tilde S \left(\|\nabla v\|_{L^{1}(M)} + \|\nabla \rho\|_{L^\infty} \|v\|_{L^1(W_\epsilon)}\right),\nonumber
\end{align}
with $\tilde S =S+S_\epsilon$. To conclude it is enough to show that
\begin{equation}\label{e:troyanov}
     \|v\|_{L^1(W_\epsilon)}\le C \|\nabla v\|_{L^{1}(M)} 
\end{equation}
for all $v\in C^\infty_c(M)$ and for some constant $C>0$ independent of $v$. To this end, take a large compact set $K\Subset M$ such that $ \mathrm{Vol}_g(K)^{\frac{n-1}{n}} \ge 2\tilde S \|\nabla \rho \|_{L^\infty}\mathrm{Vol}_g(W_\epsilon)$.
We observe first that for all $\nu \in C^\infty_c(M)$ such that $0\le \chi_ K\le \nu \le 1$ it holds 
\[
\|\nabla \nu \|_{L^1(M)}\ge \nu_0 >0.
\]
Indeed, applying \eqref{e:sobolev glued} to $\nu$ gives
\begin{align*}
\|\nabla \nu\|_{L^1(M)}&\ge \tilde S^{-1}\|\nu\|_{L^{1^\ast}(M)}- \|\nabla \rho \|_{L^\infty}\|\nu\|_{L^1(W_\epsilon)}\\
&\ge \tilde S^{-1}\mathrm{Vol}_g(K)^{\frac{n-1}{n}}- \|\nabla \rho \|_{L^\infty}\mathrm{Vol}_g(W_\epsilon)\\
&\ge \|\nabla \rho \|_{L^\infty}\mathrm{Vol}_g(W_\epsilon) =:\nu_0.
\end{align*}
Hence $\mathrm{Cap}_1(K)>0$. Take a ball $B\Subset W_\epsilon$. By Lemma \ref{lem:1cap}, $\mathrm{Cap}_1(B)=:\eta_0>0$.

From here on, one can conclude the proof as in the proof of \cite[Proposition 1]{Tr}.
%For completeness, we present the argument in our setting. 
Namely, let $v\in C^\infty_c(M)$. Without loss of generality, we can suppose that $v\ge 0$ and that \[\|v\|_{L^1(W_\epsilon)}=\mathrm{vol}_g(W_\epsilon),\]
since \eqref{e:troyanov} is invariant by rescaling. For such a rescaled $v$ it is enough to find a constant lower bound for $\|\nabla v\|_{L^1(M)}$. Thanks to the $1$-Poincar\'e inequality on $\overline{W_\epsilon}$,
\begin{equation}\label{e:poincare}
\|v - 1\|_{L^1(W_\epsilon)} \le C_W \|\nabla v \|_{L^1(M)}
    \end{equation}
for some $C_W$ which depends only on $W_\epsilon$; see \cite[Lemma 2]{Tr}. Let $\psi\in C^\infty_c(W_\epsilon)$ such that $\chi_B\le 2\psi \le 1$ and define $w\in W^{1,1}_c(M)$ by $w:=2 \max\{\psi; v\}$. Observe first that $w\ge 1$ on $B$, so that $\|\nabla w\|_{L^1(M)}\ge \eta_0.$
Define the sets
\[
A:=\{x\in W_\epsilon\,:\,\psi(x)\ge v(x)\},\qquad A':=\{x\in W_\epsilon\,:\,|v(x)-1|\ge \frac{1}{2}\}.
\]
We have that $A\subset A'$ and, by \eqref{e:poincare}, $\frac{1}{2}\mathrm{Vol}_g(A')\le C_W \|\nabla v \|_{L^1(M)}$. It holds
\[
\nabla w =\begin{cases}
    2\nabla \psi,&on\ A,\\
    2\nabla v,&on\ M\setminus A,
\end{cases}
\]
almost everywhere. In particular, a.e.,
\[
|\nabla w|\le 2|\nabla v| + 2\chi_A |\nabla \psi|,
\]
from which
\begin{align*}
    \|\nabla w\|_{L^1(M)}&\le  2\|\nabla v\|_{L^1(M)} + 2 \|\nabla \psi \|_{L^1(A)}\\
    &\le 2\|\nabla v\|_{L^1(M)} + 2 \|\nabla \psi \|_{L^\infty} \mathrm{Vol}_g(A)\\
    &\le 2\|\nabla v\|_{L^1(M)} +4C_W  \|\nabla \psi \|_{L^\infty}  \|\nabla v\|_{L^1(M)}.
\end{align*}
Thus, 
\[
\eta_0\le  \|\nabla w\|_{L^1(M)}\le  
    (2+  4C_W\|\nabla \psi \|_{L^\infty})  \|\nabla v\|_{L^1(M)}.
\]
\end{proof}

\begin{remark}\label{rmk:extension H}{\rm
 A minor variation of the above proof permits to get an analogous extension result when \eqref{e:L1Sob} is replaced by the more general inequality
\begin{equation*}%\label{e:L1Sob H}
\|h\|_{n/(n-1)}\le C \|\sqrt{|\nabla h|^2 + H^2 h^2}\|_1,
    \end{equation*}
$H$ being an assigned non-negative function. In this case \eqref{e:sobolev glued} becomes
\begin{align*}%\label{e:sobolev glued}
    \| v\|_{L^{1^\ast}(M)} &\le \|\rho v\|_{L^{1^\ast}(\Omega_\epsilon)} + \|(1-\rho)v\|_{L^{1^\ast}(M_\epsilon)}\\
    &\le S_\epsilon \|\sqrt{|\nabla(\rho v)|^2 + H^2\rho^2v^2}\|_{L^{1}(\Omega_\epsilon)} + S \|\sqrt{|\nabla((1-\rho)v)|^2+H^2(1-\rho)^2v^2}\|_{L^{1}(M_\epsilon)}\nonumber\\ 
    &\le 2\tilde S \left(\| \sqrt{|\nabla v|^2 + H^2v^2}\|_{L^{1}(M)} + \|\nabla \rho\|_{L^\infty} \|v\|_{L^1(W_\epsilon)}\right),
    %+ 2\| Hv\|_{L^1(M)},\nonumber
\end{align*}
and the proof can be concluded as above.
}\end{remark}

A straigthforward consequence of Theorem \ref{th:extension} is the following

\begin{corollary}
Suppose that $(M^n,g)$ supports the Euclidean Sobolev inequality
\begin{equation*}%\label{e:quad form}
\|v\|_{L^{\frac{n}{n-1}}(M)}\le C_M\|\nabla v\|_{L^1(M)},\quad\forall\,v\in C^\infty_c(M),
\end{equation*}
for some $C_M>0$. Let $(N^n,h)$ be such that for some $K\Subset M$ and $H\Subset N$, there is a diffeomorphism $\Phi:  M\setminus K\to N\setminus H$ with
\begin{equation}\label{e:equiv measures}
    c^{-1}\,  \Phi_\ast h \le  g \le c\, \Phi_\ast h
\end{equation}
in the sense of quadratic forms for some $c>1$. Then, there exists $C_n>0$ such that 
\begin{equation}\label{e:isop N}
    \|u\|_{L^{\frac{n}{n-1}}(N)}\le C_N\|\nabla u\|_{L^1(N)},\quad\forall\,u\in C^\infty_c(N).
    \end{equation}
\end{corollary}

\begin{proof}
    According to Theorem \ref{th:extension}, it is enough to prove the validity of \eqref{e:isop N} for $u\in C^\infty_c(N\setminus H)$. For such a $u$, set $v = u\circ \Phi \in C^\infty_c(M\setminus K)$, so that $\|v\|_{L^{\frac{n}{n-1}}(M)}\le C_M\|\nabla v\|_{L^1(M)}$. Note that $\|u\|_{L^{\frac{n}{n-1}}(N)}\le c^{n-1}\|v\|_{L^{\frac{n}{n-1}}(M)}$ and $\|\nabla v\|_{L^1(M)}\le c^{n+1}\|\nabla u\|_{L^1(N)}$, from which the thesis follows. 
\end{proof}

A special setting in which the above corollary applies is that of ALE manifolds. Following \cite{bkn}, we say that
a complete Riemannian $n$-manifold $(M^n,g)$ is asymptotically locally Euclidean
(ALE) of order $\tau > 0$, if there exists a compact subset $K \subset M$ such that
$M\setminus K$ has coordinates at infinity; namely there is $R>0$, $0<\alpha< 1$, a finite subgroup
$\Gamma\subset O(n)$ acting freely on $\mathbb R^n\setminus B_R(0)$, and a $C^\infty$-diffeomorphism $\mathcal Z:
M\setminus K \to (\mathbb R^n\setminus B_R(0))/\Gamma$ such that $\varphi = \mathcal Z^{-1} \circ \mathrm{proj}$ (where $\mathrm{proj}$ is the natural
projection of $\mathbb R^n$ to $\mathbb R^n/\Gamma$)
 satisfies \begin{eqnarray*}
      (\varphi^\ast g)_{ij}(z)=\delta_{ij} + O(|z|^{-\tau}),\qquad \partial_k(\varphi^\ast g )_{ij}(z)=O(|z|^{-\tau-1}),\\
\frac{\left|\partial_k(\varphi^\ast g )_{ij}(z)-\partial_k(\varphi^\ast g )_{ij}(w)\right|}{|z-w|^\alpha}=O(\min\{|z|,|w|\}^{-\tau-1-\alpha})\qquad\text{for}\ z,w\in \mathbb R^n\setminus B_R(0).
\end{eqnarray*}
In particular, applying the main theorem in \cite{bkn} it turns out that the isoperimetric inequality  is satisfied on $M$ if
 \[
\max |\mathrm{Sect} (x)|\le \kappa_0 r_o(x)^{-2-\epsilon}
\]
outside some compact set $K\subset M$ for some constants $\kappa_0,\epsilon>0$ and $\mathrm{AVR}>0$. 
It's worth noting that the study of isoperimetric regions in ALE manifolds or analogous classes of manifolds is an active current research field, see for instance \cite{EichmairMetzger,MondinoNardulli,Shi,ChodoshEichmairVolkmann,AntonelliNardulliPozzetta,Yu}.

Of course, the decay required in the definition of the ALE metrics is very strong. It would be interesting to have geometric (e.g. curvature) assumptions which guarantee the existence of a a diffeomorphism $\Phi:M\setminus K \to (\mathbb R^n\setminus B_R(0))/\Gamma$ satisfying \eqref{e:equiv measures}.\\

We conclude this section with an example which show that one can not deduce directly Corollary \ref{coro:nonnegoutcomp}, that is the isoperimetric inequality on manifolds whose sectional curvatures are nonnegative outside some compact set, from Theorem \ref{th:extension} and well-known isoperimetric inequalities on non-negatively curved manifolds.
\begin{example}\label{ex:cone}{\rm
Let $\mathbb C\mathbb P^2$ be the complex projective space of complex dimension 2, endowed with its Fubini-Study metric $g_{FS}$, which we can assume has sectional curvatures greater than $1$ (up to a rescaling of the metric). Let $Q=\mathrm{Con}(\mathbb C\mathbb P^2)=\mathbb R_+\times_{t} \mathbb C\mathbb P^2$ the cone over $\mathbb C\mathbb P^2$ with conical metric $dt^2 + t^2 g_{FS}$. A standard computation says that $Q$ has non-negative sectional curvature. Finally, let $M$ be the double manifold of $[1,\infty)\times_{t} \mathbb C\mathbb P^2\subset \mathrm{Con}( \mathbb C\mathbb P^2)$ endowed with a metric which coincides with $dt^2+t^2 g_{FS}$ on $P:=[2,\infty)\times_{t} \mathbb C\mathbb P^2$. Let $N$ be any complete manifold such that $P\subset N$. Since $\mathbb C\mathbb P^2$ is not null-cobordant, $N\setminus P$ can not be compact. Hence, $N$ contains a line. But then $N$ can not have non-negative Ricci curvature because of the Cheeger-Gromoll splitting theorem. Finally, note that the same argument works with $\mathbb{P}^2(\mathbb R)$ instead of $\mathrm{Con}(\mathbb{CP}^2)$. However, in that latter case we are able to deduce the validity of an isoperimetric inequality on $N$ by considering a finite (double) cover of $P$, i.e. $\mathbb R^3\setminus B_1(0)$, on which the isoperimetric inequality holds. This procedure fails with $\mathbb{CP}^2$ as this latter is simply connected.  }   
\end{example}

\section{Geometric preliminaries}\label{sec:preliminar}

Suppose that $(M,g)$ is a complete Riemannian manifold with asymptotically non-negative sectional curvature in the sense of Abresch \cite{abr}, i.e., having fixed a reference point $o\in M$,
\[
\min \sect (x) \ge -\kappa(\dist(x,o)), 
\]
for some positive monotone decreasing function $\kappa:[0,\infty)\to[0,\infty)$ satisfying 
\[
b_0(\kappa):=\int_0^{\infty} t \kappa(t)\,dt<\infty.\]
Without loss of generality, we can suppose that $\kappa$ is smooth in $[0,\infty)$ and constant in a small neighborhood of $0$.  
%$x\mapsto \kappa(\dist(x,o))$ is smooth on $M$.

Under these assumptions, Kasue \cite{ka} proved that $M$ has an asymptotic cone at infinity which is a Euclidean cone over the ideal boundary $M_\infty$ of $M$. Note that a gap in Kasue's article has been detected in \cite{drees}; however, this gap has been partially solved in \cite{mno}. For this reason, we will use in the following the approach contained in \cite{mno}, although part of the results we quote were already stated in \cite{ka}.

First, consider the rotationally symmetric Cartan-Hadamard surface $\tilde M$ with metric $    ds^2+j^2(s)d\theta^2$, with $d\theta^2$ the standard metric of $\mathbb S^1$ and $j$ a solution of 
\[
j''(s)-\kappa(s)j(s)=0,\ j(0)=0,\ j'(0)=1,\ j>0\text{ on }(0,\infty). 
\]
In particular, the Gaussian curvature of $\tilde M$ at $(s,\theta)$ is $-\kappa(s)$.
Then, consider the class $\mathcal R_o$ of rays of $M$ emanating from $o$ and $\tilde {\mathcal R}_{\tilde o}$ the class of rays emanating from $\tilde o=(0,\theta)$ in $\tilde M$.
Given $\sigma,\gamma$ in $\mathcal R_o$, $\tilde\gamma\in \tilde{\mathcal R}_{\tilde o}$ and $s,t>0$, there exists a unique ray $\tilde \sigma_{s,t} \in \tilde{\mathcal R}_{\tilde o}$ such that $\Delta\tilde o\tilde\gamma(s)\tilde\sigma_{s,t}(t)$ is a geodesic triangle of $\tilde M$ satisfying
\[
\dist_{\tilde M}(\tilde\gamma(s),\tilde\sigma_{s,t}(t))=\dist_{M}(\gamma(s),\sigma(t)). 
\]
Moreover, the angle $\angle \tilde\gamma(s)\tilde o\tilde\sigma_{s,t}(t)$ is monotone non-increasing as $t$ and $s$ increase; see \cite[Theorem 1.3]{IMS} and \cite[Theorem 1.1]{mno}. In particular, $\tilde \sigma_{s,t}$ converges to a ray $\tilde \sigma_ {\infty}\in \tilde{\mathcal R}_{\tilde o}$ as $s,t\to\infty$, in the sense of the uniform convergence on bounded sets.
Thus, we can introduce the angle distance $\angle_\infty$ on $\mathcal R_o$ by 
\begin{equation}\label{e:MNO1}
\angle_\infty(\gamma,\sigma):= \angle_\infty^{\tilde M}(\tilde\gamma,\tilde\sigma_\infty),
    \end{equation}
where the angle distance at infinity on $\tilde M$ is given by the limit as $s,t\to\infty$ of the angle at the origin in the flat comparison triangle, i.e.
\begin{equation}\label{e:MNO2}
\angle_\infty^{\tilde M}(\tilde\gamma,\tilde\sigma_\infty):=\lim_{t,s\to\infty}\arccos\frac{t^2+s^2-\dist_{\tilde M}(\tilde\gamma(s),\tilde \sigma_\infty(t))^2}{2ts},
\end{equation}
and is well defined as the argument of the limit is monotone non-decreasing because of the non-positiveness of the curvature of  $\tilde M$; see \cite[(1.1)]{mno}.
Finally, we define the ideal boundary $M_\infty$ of $M$ as the quotient of $\mathcal R_o$ via identification of equivalent rays, that is rays whose distance $\angle_\infty$ is zero. The space $M_\infty$ has a metric structure induced by $\angle_\infty$ through the quotientation process, which we still denote by $\angle_\infty$ with the usual abuse of notation. Moreover, a second metric structure on $M_\infty$ is given by the intrinsic metric $\angle_{Tits}$ induced by $\angle_\infty$. 

In the following lemma we collect a number of useful facts about the ideal boundary of asymptotically nonnegatively curved manifolds.

\begin{lemma}\label{lem:mno}Suppose that $(M,g)$ is a complete Riemannian manifold with asymptotically non-negative sectional curvature in the sense of Abresch. Then
    \begin{enumerate}
        \item For all $[\sigma],[\gamma]\in M_\infty$, $\angle_\infty([\sigma],[\gamma])=\min \{\angle_{Tits}([\sigma],[\gamma]);\pi\}$.
        \item For all $\sigma,\gamma\in \mathcal R_o$ and all $t>0$, 
        \[
\angle_\infty([\sigma],[\gamma])=\lim_{t\to\infty}2\arcsin\left( \frac{\dist(\sigma(t),\gamma(t))}{2t}\right).
\]
\item Let $M_{\infty,0}$ be a connected component of
the Tits ideal boundary of $M$ and let us denote by $\mathrm{Cone}(M_{\infty,0})$ the Euclidean cone over $M_{\infty,0}$ with vertex $o^{*}$. Then $\mathrm{Cone}(M_{\infty,0})\setminus\{o^*\}$ is an Alexandrov space with curvature bounded below by $0$, and $\mathrm{dim}_{\mathcal H} (M_{\infty,0})$ is an integer not greater than $\mathrm{dim}(M)-1$. More precisely,
the following hold:
\begin{enumerate}
    \item If $\mathrm{dim}_{\mathcal H} M_{\infty,0} = 1$ and if the diameter of $M_{\infty,0}$ is not greater than $\pi$, then $\mathrm{Cone}(M_{\infty,0})$ itself
and $M_{\infty,0}$ are Alexandrov spaces with curvature bounded below by $0$ and $1$, respectively.
\item If $\mathrm{dim}_{\mathcal H} M_{\infty,0}\ge 2$, then $\mathrm{Cone}(M_{\infty,0})$ and $M_{\infty,0}$ are Alexandrov spaces with curvature
bounded below by $0$ and $1$, respectively. In particular, $\diam M_{\infty,0}\le \pi$.
\end{enumerate}
\item $\mathrm{dim}_{\mathcal H} M_{\infty}=\mathrm{dim}(M)-1$ if and only if $\textrm{AVR}>0$.
    \end{enumerate}
\end{lemma}
\begin{proof} We refer to \cite{mno} for the proofs. More precisely
    \begin{enumerate}
        \item See \cite[Theorem 2.5]{mno} or \cite[Proposition 2.4]{ka}.
        \item See \cite[Proposition 2.1]{mno}.
        \item See \cite[Corollary 0.6]{mno}.
        \item See \cite[Corollary 0.7]{mno}.
    \end{enumerate}
\end{proof}

We note for later purposes that each connected component $M_{\infty,0}$ of $M_\infty$ is compact. If $n\ge 3$, then $M_{\infty,0}$ has curvature not less than $1$ and its diameter is not greater than $\pi$, because of the Bonnet-Myers theorem for Alexandrov spaces. Fix $\alpha\in(0,\frac{\pi}{3})$. Then, by the Bishop-Gromov inequality, for every $\sigma\in \mathcal R_o$,
\begin{equation*}%\label{e:vol}
\mathcal H_{n-1}(B_{\alpha}([\sigma]))\ge
\frac{\mathcal H_{n-1}(B_\alpha^{\mathbb S^{n-1}})}{\mathcal H_{n-1}(\mathbb S^{n-1})}
\mathcal H_{n-1}(M_{\infty,\sigma}),    \end{equation*}
where  $M_{\infty,\sigma}$ is the connected component of $M_\infty$ containing $[\sigma]$.
    Instead, if $n=2$, $M_{\infty,\sigma}$ is a circle and 
\begin{equation*}%\label{e:vol2}
\mathcal H_{1}(B_{\alpha}([\sigma]))\ge \min\{\mathcal H_{1}(M_{\infty,\sigma});\alpha\}.
    \end{equation*}
Let $E_\sigma$ be the end of $M$ containing (part of) $\sigma$, and let $\mathrm{Cone}_\sigma$ be the metric cone over $M_{\infty,\sigma}$. According to (the proof of) \cite[Corollary 0.7]{mno}, 
\[
\lim_{r\to\infty}\frac{\vol(B_r(o)\cap E_\sigma)}{r^n} = \mathcal H_n(B_1^{\mathrm{Cone}_\sigma}) = \frac{1}{n}\mathcal H_{n-1}(M_{\infty,\sigma}).\]
Recapping, for every $n\ge2$ and $\sigma\in\mathcal R_o$, it holds
\begin{equation}\label{e:vol}
\mathcal H_{n-1}(B_{\alpha}([\sigma]))\ge
\mu(\alpha,\theta,n),\end{equation}
for some constant $\mu$, where \[\theta:=\inf_{E_j} \limsup_{r\to\infty}\frac{\vol(B_r(o)\cap E_j)}{r^n}\] is positive, and thus $\mu(\alpha,\theta,n)>0$ whenever all the ends of $M$ have positive $\mathrm{AVR}$.

We also need the following

\begin{lemma}\label{lem:Lip}
For every $\eta\in(0,1)$ there exists a radius $\tilde T=\tilde T(\eta,b_0(\kappa))$ such that
\[
\frac{1}{t}\dist(\sigma(t),\gamma(t))\ge  2\sin\left(\frac{1-\eta}{2}\angle_\infty([\sigma],[\gamma])\right),
\]
for all $\sigma,\gamma\in \mathcal R_o$ and all $t>\tilde T$, as soon as $\dist(\sigma(t),\gamma(t))<\tilde T$. In particular, there exists $C_\eta>0$ such that if $t>\tilde T$ and $\dist(\sigma(t),\gamma(t))<C_\eta\,\tilde T$ then 
\[
\frac{1}{t}\dist(\sigma(t),\gamma(t))\ge  (1-\eta)^2\angle_\infty([\sigma],[\gamma]).
\]
\end{lemma}

\begin{proof}
Let $\eta\in(0,1)$. Choose $\eta_{0}\in(0,1)$ such that
\[
\ (1-\eta_{0})\sin u\geq \sin ((1-\eta)u),\quad \forall\,u\in[0,\pi/2].
\]
Let $\sigma,\gamma\in \mathcal R_o$ and $s,t>0$. Consider the geodesic triangle $\triangle \tilde\gamma(t)\tilde o \tilde\sigma_{t,t} (t)$ in $\tilde M$ given by  \cite[Theorem 1.3]{IMS}, with $\tilde \sigma_{t,t},\tilde\gamma\in\tilde{\mathcal R}_{\tilde o}$ and $\dist_{\tilde M}(\tilde\gamma(t),\tilde\sigma_{t,t} (t))=\dist_{M}(\gamma(t),\sigma(t))$.
Since $\angle \tilde\gamma(t)\tilde o\tilde\sigma_{t,t}(t)$ is monotone non-increasing in $t$, it holds that $\dist_{\tilde M}(\tilde\gamma(t),\tilde\sigma_{t,t} (t))\ge \dist_{\tilde M}(\tilde\gamma(t),\tilde\sigma_{\infty} (t))$.  Indeed, in polar coordinates, the map $(t,\theta)\mapsto(t,\alpha\theta)$, $0<\alpha<1$ is $1$-Lipschitz on any sector and hence, in the rotationally symmetric model, the distance between two points on the same geodesic sphere is nondecreasing with respect to the angle.

We claim that for $\tilde T$ large enough, for every $s>t>\tilde T$, if $\dist_{\tilde M}(\tilde\gamma(t),\tilde\sigma_{\infty} (t))<\tilde T$ then
\begin{equation}\label{e:lip}
\frac{\dist_{\tilde M}(\tilde\gamma(t),\tilde\sigma_{\infty} (t))}{t}\ge (1-\eta_{0})\frac{\dist_{\tilde M}(\tilde\gamma(s),\tilde\sigma_{\infty} (s))}{s}
\end{equation}
To this end, note that by L'Hôpital rule $j(t)/t$ tends to a constant $j'(\infty)=\lim_{t\to\infty}j'(t)$ as $t\to\infty$; see e.g. \cite[Lemma 2.1]{abr}. 
Hence, for $\tilde T$ large enough $j(t)>(1-\eta_{0})\frac{t}{s}j(s)$ for all $s>t>\frac{\tilde T}{2}$. In particular the maps of $\tilde M$ given in polar coordinates by $(\rho,\theta)\mapsto (\frac{s}{t}\rho,\theta)$ is $\frac{s}{t(1-\eta_{0}})$-Lipschitz on $    \tilde M^{\tilde T/2}:=\tilde M\setminus B_{\tilde T/2}(\tilde o)$. Since $\tilde M^{\tilde T/2}$ contains the minimizing geodesics of $\tilde M$ realizing $\dist_{\tilde M}(\tilde\gamma(t),\tilde\sigma_{\infty} (t))$, \eqref{e:lip} is proved. 

Letting $s\to\infty$ and recalling \eqref{e:MNO1} and \eqref{e:MNO2} gives the first statement. The second statement is a straightforward consequence of the first one. Indeed,  if $C_\eta= 2\sqrt{\eta}(1-\eta)$, then $2\arcsin(\frac{x}{2})\le \frac{x}{1-\eta}$ for all $x\in[0,C_\eta]$, whence the thesis.
\end{proof}

Here and in the following we let $r_o(x):=\dist(o,x)$. Moreover, given $u\in S_oM$, we denote by $\gamma_u(t)=\exp_o(tu)$  the unit speed geodesic leaving $o$ with velocity $u$, defined on the maximal set $[0,\ell)$ where $\gamma_u$ is minimizing.

\begin{lemma}\label{lem:topo2}
Let $a,\epsilon\in(0,1)$. Let $\gamma_u,\gamma_v\in \mathcal R_o$ two rays with $\dot\gamma_u(0)=u,\dot\gamma_v(0)=v\in S_oM$. Then
\begin{align}\label{e:topo}
    \dist(\gamma_{u}(t),\gamma_{v}(t)) &\le  \left(1+\epsilon-(1-\epsilon)\sqrt{1-a^2}\right) t,
\end{align}
for all $t>0$ as soon as $\angle_o(u,v)\le \arccos\sqrt{1-a^2\epsilon^2\beta^2}=:\alpha_{0}$, where $\beta$ is the positive finite constant depending on $\kappa$ introduced in \cite[p.659]{abr}.
In particular, \[
\angle_\infty([\gamma_u],[\gamma_v])\le \left(1+\epsilon-(1-\epsilon)\sqrt{1-a^2}\right)
\]
\end{lemma}
\begin{proof}
    This is a direct consequence of \cite[Proposition 1]{abr} and of the triangle inequality.
\end{proof}
As a consequence we have that the map $u\mapsto [\gamma_u]$ is continuous when restricted to the set $\mathcal S_\infty:=\{u\in S_oM\ :\ \gamma_u\in \mathcal R_o\}$ where it is well-defined. 
\begin{lemma}\label{lem:unif lim}
    The limit in Lemma \ref{lem:mno} \textit{(2)} is uniform: for every $\bar \epsilon>0$ there exists $T>0$ such that for all $t\ge T$ and all $\sigma_1,\sigma_2\in\mathcal R_o$
    \[
\left|\frac{\dist(\sigma_{1}(t),\sigma_{2}(t))}{t} - 2\sin\left(\frac{\angle_\infty([\sigma_{1}],[\sigma_{2}])}2\right)\right|<\bar \epsilon.
\]
\end{lemma}
\begin{proof}
Suppose by contradiction that for some $\bar \epsilon>0$ there exists a sequence $\{t_k\}_{k=1}^\infty$ with $\lim_kt_k=\infty$ and a corresponding sequence $(\sigma_{1,k},\sigma_{2,k})\in \mathcal R_o\times\mathcal R_o$ such that 
\[
\left|\frac{\dist(\sigma_{1,k}(t_k),\sigma_{2,k}(t_k))}{t_k} - 2\sin\left(\frac{\angle_\infty([\sigma_{1,k}],[\sigma_{2,k}])}2\right)\right|>\bar \epsilon.
\]
Up to extract a subsequence, by compactness we can assume that $\dot\sigma_{i,k}(0)$ converges to $u_{i,\infty}$ in $S_oM$ as $k\to\infty$, $i=1,2$. Accordingly, by Lemma \ref{lem:topo2} and the subsequent observation we have that \[\angle_\infty([\sigma_{1,k}],[\sigma_{2,k}])\to \angle_\infty([\gamma_{u_{1,\infty}}],[\gamma_{u_{2,\infty}}]), \qquad k\to\infty.
\]
However, for $k$ large enough
\[
\frac{\dist(\sigma_{i,k}(t_k),\gamma_{u_{i,\infty}}(t_k))}{t_k}\le \frac{\bar \epsilon} 3,\qquad i=1,2,
\]
so that
\[
\left|\frac{\dist(\gamma_{u_{1,\infty}}(t_k),\gamma_{u_{2,\infty}}(t_k))}{t_k} - 2\sin\left(\frac{\angle_\infty([\sigma_{1,k}],[\sigma_{2,k}])}2\right)\right|>\frac{\bar \epsilon}{3}.
\]
This latter gives a contradiction in the limit as $k\to\infty$.
\end{proof}

The following general fact holds on any complete non-compact Riemannian manifold. Here by $\angle_o (\gamma,\sigma)$ we denote the angle at $o$ between $\gamma$ and $\sigma$.
\begin{lemma}\label{lem:density radial}
For any $\zeta>0$ there exists a positive radius $\rho_0$ such that for any $x\in M\setminus B_{\rho_0}(o)$ and any minimizing geodesic $\gamma$ from $o$ to $x$ there exists a ray $\sigma\in\mathcal R_o$ with $\angle_o (\gamma,\sigma)<\zeta$.
\end{lemma}

\begin{proof}
Let $\mathcal U:=\{U_{k}\}_{k=1}^K$ be a finite open covering of $S_oM$, the space of unitary tangent vectors at $o$, such that, for all $k$, $\diam U_{k}<\zeta$, where the diameter is computed with respect to the spherical distance on $S_oM\asymp \mathbb S^{n-1}$. Since $M$ is non-compact there is at least a $\bar k$ such that $\overline{U_{\bar{k}}}$ contains a vector which generates a ray. For any such $\bar k$, set $R_{\bar{k}}=0$. For any $\overline{U_k}$ generating no rays, let 
\[R_k:=\sup_{w\in \overline{U_k}}\sup\{t>0\,:\,r_o(\exp_o(tw))=t\}.\]
By continuity of the exponential map, $R_k$ is finite for all $k$. Hence the assertion of the lemma holds true by choosing $\rho_0=\max_k R_k$.
\end{proof}

\section{Asymptotically non-negative curvature
}\label{sec:asymptotic}

This section is devoted to the proof of 
Theorem \ref{thm:asymp nonneg}.
%     Let $M$ be 
%  a complete smooth  $n$-dimensional Riemannian manifold and let $o\in M$. Suppose that
%  \[
% \min \mathrm{Sect} (x)\ge -\kappa_0 r_o(x)^{-n-\epsilon}
% \]
% outside some compact set $K\subset M$ for some constants $\kappa_0,\epsilon>0$. 
%  Suppose that   $\textrm{AVR}>0$ on any end $E_j$ of $M$, i.e. 
%     \[
%     \limsup _{r\to\infty}\frac{\vol(B_r(o)\cap E_j)}{r^n}>0.
%     \] 
%     Then, there exists a constant $C=C(M)>0$ such that for any compact domain $\Omega\subset M$ with smooth boundary and any $h$ a positive smooth function on $\Omega$ it holds
%     \begin{equation}\label{e:isoper asymp}
%             \int_\Omega |\nabla h| + \int_{\partial \Omega} h \ge C \left(\int_\Omega h^{\frac{n}{n-1}}\right)^{\frac{n-1}{n}}.
%         \end{equation}
%\end{theorem}

First, we observe the smallness (in a spectral sense) of the negative part of the curvature of $M$. This will be used to conformally deform the metric of $M$ into a new metric of non-negative weighted curvature in the exterior of a compact set.

\begin{proposition}\label{prop:supersol}
    Let $(M^n,g)$ be a complete non-compact Riemannian manifold whose sectional curvatures satisfy
    \[
    \min \mathrm{Sect} (x) \ge -\kappa( r_o(x))
    \]
    for some $o\in M$ and some positive non-increasing function $\kappa$ such that
    \[ \kappa(t) \le \kappa_0 t^{-(n+\epsilon)},\qquad \forall\,t\ge t_0\]
    for some $\kappa_0,t_0,\epsilon>0$. Then, for any $\lambda>0$ there exists a function $\varphi\in C^\infty(M)$ which satisfies
    \begin{equation}\label{e:critical}
    \begin{cases}
        \Delta \varphi (x) + \lambda (n-1) \kappa(r_o(x)) \varphi(x) \le 0, \\
    |\varphi(x)-1+r_o(x)^{-(n-2+\epsilon)/2}|\le r_o(x)^{-2}
    \end{cases}
    \end{equation}
     on $M\setminus B_{t_1}(o)$
    for some $t_1>0$. 
\end{proposition}

\begin{proof}
    First, take $A>0,0<B<D'<D$ constants and note that 
    $y(t)=1-t^{-2D'}$ solves    
    the ordinary differential inequality
    \begin{equation}\label{e:ODE}
    y''(t)+ (1+2B)\frac{y(t)}{t}+A^2\frac{y(t)}{t^{2+2D}}\le - C t^{-2-2D'}<0
    \end{equation}
   on $[T_0,\infty)$ with $C=2D'(D'-B)>0$ for some $T_0$ large enough depending on $A,\, B,\, D,\, D'$. Moreover $y'>0$. Then, define
   $v$ on $M\setminus B_{T_0}(o)$ by $v(x):=y(r_o(x))$. Observe that, as $y$ is increasing and $\mathrm{Ric}(x)\ge -(n-1)\kappa(r_0(x))$, by the Laplacian comparison \cite[Theorem 2.4]{prs}
   \[
   \Delta v (x) \le  y''(r_o(x))+(n-1)\frac{\sigma'(r_o(x))}{\sigma(r_o(x))}y'(r_o(x))
   \]
   weakly, where
   $\sigma:[0,+\infty)\to [0,+\infty)$ is the function solving
   \[
   \begin{cases}
       \sigma''(t)-\kappa(t) \sigma(t) = 0,\\
       \sigma(0)=0,\quad 
       \sigma'(0)=1.
   \end{cases}
   \]
   In particular, $t\sigma'(t)/\sigma(t)\to 1$ as $t\to\infty$. Indeed, 
\[
  \lim_{t\to\infty}\frac{(t\sigma'(t))'}{\sigma'(t)}=1 + \lim_{t\to\infty}\frac{t\kappa(t)\sigma(t)}{\sigma'(t)}=1
  \]
  since $\sigma(t)\le e^{b_0(\kappa)}t$,  $\sigma'(t)\ge 1$ and $\kappa(t)=o(t^{-2})$. The claim follows from L'H\^opital's rule.

   Now, choose $(n-2)/2<D<(n-2+\epsilon)/2$, $B\in((n-2)/2,D)$, $A^2=(n-1)\lambda\kappa_0$. Thus, there exists $t_1\gg 1$ such that 
   \begin{align}
          \Delta v(x) + \lambda (n-1) \kappa(r_o(x)) v(x) &\le y''(r_o(x))+(n-1)\frac{\sigma'(r_o(x))}{\sigma(r_o(x))}y'(r_o(x)) + \lambda (n-1) \kappa(r_o(x)) y(r_o(x))\label{e:weak supersol}\\
          &\le y''(r_o(x))+ (1+2B)\frac{y'(r_o(x))}{r_o(x)} + A^2 \frac{y(r_o(x))}{r_o(x)^{2+2D}}\nonumber\\
          &\le -Cr_o(x)^{-2-2D'}=:-\psi(x),\nonumber
      \end{align}
      weakly on $M\setminus B_{t_1}(o)$. Let $\tilde{\kappa}\doteq\kappa\circ r_o$. The existence of a supersolution for the operator $L:=\Delta + \lambda (n-1) \tilde \kappa$ guarantees that the bottom of the spectrum of $-L$ on $M\setminus B_{t_1}(o)$ is positive. Hence, there exists a function $0<\alpha$ which solves $L\alpha =0$ on $M\setminus B_{t_1}(o)$, \cite{FiCoSh}. Using a well-known trick introduced by M.H. Protter and H.F. Weinberger, \cite{PW}, we deduce that \eqref{e:weak supersol} is equivalent to
      \[
      \Delta\left(\frac{v}{\alpha}\right)+2g\left( \frac{\nabla\alpha}{\alpha},\nabla\left(\frac{v}{\alpha}\right)\right) \le -\frac{\psi}{\alpha}<0
      \]
      weakly on $M\setminus B_{t_1}(o)$; see \cite[Lemma 2.3]{PVV}.

      To get rid of the drift, we can use a second well-known trick. Namely, set $\rho=\alpha^{2/n}$ and consider the warped product $N=M\times_{\rho}\mathbb T^1$  endowed with the warped metric $g_\rho:=g+ \rho^2 g_{\mathbb T^1}$. Here, $\mathbb T^1=\mathbb R / \mathbb Z$. An explicit computation shows that \[
      \Delta_N  w=\Delta w  +n g_{\rho}( \nabla \log \rho,\nabla w) + \frac{1}{\rho^2}\partial^2_{\theta\theta} w 
   \]
   where $\theta$ is the coordinate of $\mathbb T^1$ and $\Delta_N$ is the (unweighted) Laplace-Beltrami operator on $(N,g_\rho)$. Choosing $w(x,\theta)=-v(x)/\alpha(x)$, we obtain in particular that $\Delta_N w\ge \psi/\alpha >0$ weakly on $N$.
An 
       application of \cite[Corollary 1]{GW} guarantees the existence of a $C^\infty$ strictly subharmonic (with respect to $g_\rho$) function $z$ on $(M\setminus B_{t_1}(o)\times_\rho \mathbb T^1,g_\rho)$ such that 
      \[
\left| z(x,\theta) + \frac{v(x)}{\alpha(x)} \right|  =      | z(x,\theta) - w(x,\theta) |     <  \frac{r_o(x)^{-2}}{\alpha(x)}
      \]
      and
 \[
 \Delta z(x,\theta) + n g_{\rho}( \nabla \log \rho, \nabla z (x,\theta)) \ge - \frac{1}{\rho^2}\partial^2_{\theta\theta} z.
 \]     
Finally, define $\varphi (x)= -\alpha (x)\int_{\mathbb T^1} z(x,\theta)\,d\theta$. Clearly, $|\varphi (x) - v(x) |  <  r_o(x)^{-2}$. Moreover,
\begin{align*}
    \Delta \frac{\varphi}{\alpha} + n g_{\rho}( \nabla \log \rho, \nabla \frac{\varphi}{\alpha})
    &=-\Delta \int_{\mathbb T^1} z(x,\theta)\,d\theta - n g_{\rho}( \nabla \log \rho, \nabla \int_{\mathbb T^1} z(x,\theta)\,d\theta)\\
&= -\int_{\mathbb T^1} \left\{\Delta z(x,\theta) + n g_{\rho}( \nabla \log \rho, \nabla  z(x,\theta))   \right\} \,d\theta\\
&\le   \frac{1}{\rho^2(x)} \int_{\mathbb T^1} \partial^2_{\theta\theta} z(x,\theta)=0,\end{align*}
so that 
\begin{align*}
\Delta \varphi (x) + \lambda (n-1) \kappa(r_o(x)) \varphi (x)
\le
\Delta \varphi (x) + \lambda (n-1) \tilde \kappa (x) \varphi (x) \le 0.
      \end{align*}
 \end{proof}
\begin{remark}{\rm
    Sharper estimates can be obtained considering the exact solution to 
    \begin{equation*}%\label{e:ODE}
    y''(t)+ (1+2B)\frac{y(t)}{t}+A^2\frac{y(t)}{t^{2+2D}}=0
    \end{equation*}
    converging to $1$ at infinity, i.e. 
    \begin{equation}\label{e:Bessel}
        y(t)=(2D/A)^{B/D}\Gamma(1-B/D)t^{-B}J_{-B/D}(At^{-D}/D)
    \end{equation} where $J_{-B/D}$ is the Bessel function of the first kind of index $-B/D$. Note however that this choice permits to improve the constants, but not the range of admissible coefficients. When $B>D$ (situation which corresponds to a decay of $\kappa$ with a power exponent between $2$ and $n$), then $y$ as in \eqref{e:Bessel} is still a solution, however it turns out to be decreasing so that Laplacian comparison cannot be applied. We do not know if and when one can guarantee that a solution to \eqref{e:critical} exists in this case.
}\end{remark}

We can now prove Theorem \ref{thm:asymp nonneg}.

\begin{proof}
Let $r_o(x):=\dist(o,x)$. Fix a small constant $\delta>0$ so that
$2 \sin \frac\alpha 2 < 1 - 2\delta$, where $\alpha\in(0,\frac\pi 3)
$ is the constant introduced in Section \ref{sec:preliminar}.

Choose $\epsilon,a\in (0,1)$ small enough so that 
\[
(1+\epsilon-(1-\epsilon)\sqrt{1-a^2})<\delta/2,
\]
and set $\alpha_0=\arccos\sqrt{1-a^2\epsilon^2\beta^2}$ as in Lemma \ref{lem:topo2}.
Let $\rho_0$ be the radius from Lemma \ref{lem:density radial} applied with $\zeta:=\alpha_0$.

Fix $\lambda=\max\{n-2;1\}$ and let $t_1$ and $\varphi$ be the radius and the function given by Proposition \ref{prop:supersol}. Recall that $\varphi\in C^\infty(M)$ solves
\[
\Delta \varphi + \lambda(n-1)(\kappa\circ r_0) \varphi \le 0
\]
on $M\setminus B_{t_1}(o)$ and \[
1-r_o(x)^{-(n-2+\epsilon)/2} - r_o(x)^{-2} \le \varphi(x) \le 1-r_o(x)^{-(n-2+\epsilon)/2} + r_o(x)^{-2}
\]
therein. In particular, up to increase $t_1$ and using that  \[b^{1/\lambda}-a^{1/\lambda}=\int_a^b \frac {s^{-1+{1/\lambda}}}{\lambda }ds\le \frac{(b-a)}{\lambda}\left(\min_{[a,b]}s\right)^{-1+1/\lambda},\qquad b>a>0,\] we have \begin{align}\label{e:est vp lambda}
\left|\varphi(x)^{1/\lambda}-(1-r_o(x)^{-\frac{n-2+\epsilon}{2}})^{1/\lambda}\right| &\le \left| \varphi(x) - 1 +r_o(x)^{-\frac{n-2+\epsilon}{2}}\right|\frac{\left(1-r_o(x)^{-\frac{n-2+\epsilon}{2}}\right)^{-1+1/\lambda}}{\lambda}\\
&\le \frac{1}{2\lambda} r_o(x)^{-2}\nonumber
\end{align}
on $M\setminus B_{t_1}(o)$. Moreover, without loss of generality, we can suppose that $t_1\ge  \max\{\rho_0,T\}$, where $T$ is the radius in Lemma \ref{lem:unif lim} applied with $\bar\epsilon = \delta$.

Finally, let $T_0\ge t_1$ large enough so that
\begin{equation}\label{e:T0}
   \left(1-\frac{\delta}{2}\right) t+ \frac{2}{\lambda}\int_{t_1}^\infty s^{-2}\,ds \le  \int_{t_1}^{t}(1-s^{-(n-2+\epsilon)/2})^{1/\lambda}\,ds,\qquad \forall\,t\ge T_0.
\end{equation}

Following \cite{CMT,IRV}, 
set $f:=\log \varphi / \lambda$ and consider on $M$ the conformal metric $\tilde g=e^{2f}g$. Moreover, for the ease of notation define $f_n:=(n-2)f$.

We are going to prove that the isoperimetric inequality holds on $M\setminus B_{T_0}(o)$. Thanks to Theorem \ref{th:extension} this will permit to conclude. 

In what follows, the superscript $\tilde{\ }$ will denote objects which are defined in terms of the metric $\tilde g$.

 Let $\Omega\subset M\setminus B_{T_0}(o)$ be a compact domain with smooth boundary $\partial \Omega$ and let $h$ be a positive smooth function on $\Omega$. Preliminarily, remark that we can assume without loss of generality that $\Omega$ is connected, as the unconnected case follows easily from the connected one; see \cite[p.6]{brendle-Rn}. 
We assume by scaling that
\begin{equation}\label{e:scalingconf}   
\int_\Omega |\widetilde{\nabla h}|_{\tilde g}\,e^{-f_n}\,d\widetilde{\mathrm{vol}} + \int_{\partial \Omega} h \,e^{-f_n}\,d\widetilde{\mathrm{vol}_{n-1}} = n \int_\Omega h^{\frac{n}{n-1}}\,e^{-f_n}\,d\widetilde{\mathrm{vol}},
\end{equation}
and we consider the Neumann problem
\[
\begin{cases}
e^{f_n}\widetilde{\mathrm{div}}(e^{-{f_n}}h\widetilde{\nabla} u)= nh^{\frac{n}{n-1}}-|\widetilde{\nabla} h|_{\tilde g},&\textrm{in }\Omega\\
\partial_{\tilde\nu} u = 1,& \text{on }\partial\Omega.
\end{cases}
\]
The solution is continuous on $\bar\Omega$ and unique up to an additive constant. Hence, there exists a point $x_\ast\in\Omega$ where $u$ attains its minimum. Let $\sigma_\ast(t)=\exp_o(t w_\ast)$ be a minimizing geodesic from $o$ to $x_\ast=:\sigma_\ast(t_\ast)$. Note that $w_\ast\in\mathcal S_{\rho_0}$, so that $\angle_o (\sigma_\ast,\gamma_\ast)<\alpha_0/2$ for some $\gamma_\ast\in\mathcal R_o$. 
Set $V=B_{\alpha}([\gamma_\ast])\subset M_\infty$.

\begin{lemma}\label{lem:distances}
    For any ray $\sigma\in \mathcal{R}_0$ with $[\sigma]\in V$ it holds
    \[
     \widetilde{\dist}(\sigma(t),B_{t_1}(o))\ge  \widetilde{\dist}(\sigma(t),x_\ast) ,\qquad \forall t>t_\ast.
    \]
\end{lemma}
\begin{proof}[Proof (of Lemma \ref{lem:distances})]
    First, by \eqref{e:topo} and the choice of $\alpha_0$ we have that
    \[
    \dist(x_\ast,\gamma_{\ast}(t_\ast)) \le \frac{\delta}{2} t_\ast .
    \]
    Since $t_\ast>T$, from Lemma \ref{lem:unif lim} applied with $\bar\epsilon=\delta$ and the choice of $T$, we also have
    \[
    \dist(\gamma_{\ast}(t_\ast),\sigma(t_\ast))<\left[2\sin\left(\frac{\alpha}2\right)+\delta\right]t_\ast< \left(1-\delta \right)t_\ast,\quad \forall \sigma\in V.
    \]
    
Hence, by the triangle inequality, 
\[   \widetilde{\dist}(x_\ast,\sigma(t_\ast))\le \dist(x_\ast,\sigma(t_\ast)) \le \left(1-\frac{\delta}{2}\right)t_\ast,\]
 where we also used that $e^{f}\le 1$, so that $\widetilde{\dist}\le\dist$. Now, using \eqref{e:est vp lambda}, for any $t>t_\ast$,
 \[
 \widetilde{\dist}(\sigma(t_\ast),\sigma(t))\le \int_{t_\ast}^t e^{f(\sigma(s))}\,ds\le \int_{t_\ast}^t \left((1-s^{-(n-2+\epsilon)/2})^{1/\lambda} + \frac{s^{-2}}{2\lambda}\right)\,ds.
 \]
 Hence, by the triangle inequality, for any $t>t_\ast$ we obtain
 \begin{align*}
 \widetilde{\dist}(x_\ast,\sigma(t))&\leq \widetilde{\dist}(x_\ast,\sigma(t_\ast))+ \widetilde{\dist}(\sigma(t_\ast),\sigma(t))\\
 &\leq \left(1-\frac{\delta}{2}\right)t_\ast+ \int_{t_\ast}^t \left((1-s^{-(n-2+\epsilon)/2})^{1/\lambda} + \frac{s^{-2}}{2\lambda}\right)\,ds.
 \end{align*}
 On the other hand, let $q\in\partial B_{t_1}(o)$ ($B_{t_1}(o)$ being here the geodesic ball with respect to the metric $g$) be a point that realizes the distance for the metric $\tilde g$ from $\sigma(t)$ to $\partial B_{t_1}(o)$ and let 
 $\gamma:[0,\ell]\to M\setminus B_{t_1}(o)$ be a minimizing geodesic for the metric $\tilde g$ joining $\gamma(0)=q$ to $\gamma(\ell)=\sigma(t)$, parametrized so that $g(\dot\gamma,\dot\gamma)=1$. Note that
 \[
 \ell=\mathrm{Length}_g(\gamma)\geq \mathrm{dist}(\sigma(t),q)\geq r_o(\sigma(t))-r_o(q)=t-t_1.
 \]
 Then,
 \begin{align*}
    \widetilde{\dist}(\sigma(t),B_{t_1}(o))&=\mathrm{Length}_{\tilde g}(\gamma)\ge \int_0^\ell e^{f(\gamma(s))}\,ds\\
    &\ge\int_{\ell-(t-t_1)}^\ell e^{f(\gamma(s))}\,ds\\
    &=\int_{0}^{t-t_1} e^{f(\gamma(\ell-q))}\,dq\\
    &\ge\int_{0}^{t-t_1}\left((1-(t-q)^{-(n-2+\epsilon/2)})^{1/\lambda} - \frac{(t-q)^{-2}}{{2\lambda}}\right)\,dq\\
& = \int_{t_1}^{t}\left((1-s^{-(n-2+\epsilon)/2})^{1/\lambda} - \frac{s^{-2}}{2\lambda}\right)\,ds,
  \end{align*}
where we used that $r_o(\gamma(\ell-q))\ge r_0(\gamma(\ell))-\dist(\gamma(\ell-q),\gamma(\ell))\ge t-q$ and \eqref{e:est vp lambda} once again.
Accordingly,
 \[
 \widetilde{\dist}(\sigma(t),B_{t_1}(o))\ge  \widetilde{\dist}(x_\ast,\sigma(t)) 
 \]
 provided that
 \begin{align*}
     \left(1-\frac{\delta}{2}\right)t_\ast+ \int_{t_\ast}^t \left((1-s^{-(n-2+\epsilon)/2})^{1/\lambda} + \frac{s^{-2}}{2\lambda}\right)\,ds &\le \int_{t_1}^{t}\left((1-s^{-(n-2+\epsilon/2)})^{1/\lambda} - \frac{s^{-2}}{2\lambda}\right)\,ds 
  \end{align*}
 which in turn is implied by
 \[
\left(1-\frac{\delta}{2}\right)t_\ast+ \frac{2}{\lambda}\int_{t_1}^\infty s^{-2}\,ds \le \int_{t_1}^{t_\ast}(1-s^{-(n-2+\epsilon)/2})^{1/\lambda}\,ds.
 \]
 This latter is guaranteed by \eqref{e:T0}, as $t_\ast \ge T_0$.
\end{proof}

As in \cite{brendle}, define
 \[
 U=\{x\in\Omega\setminus\partial\Omega\,:\,|\widetilde{\nabla} u(x)|_{\tilde g}<1\}
 \]
 and, for $r>0$,
 \[
A_r=\{x\in U\,:\,\forall y\in U,\ ru(y)+\frac{1}{2}\widetilde{\dist}(y,\widetilde{\exp}_x(r\widetilde{\nabla} u(x))^2 \ge ru(x) +\frac{1}{2}r^2|\widetilde{\nabla} u(x)|_{\tilde g}^2\}. 
 \]
Define the $C^{1,\alpha}$
transport map $\widetilde{\Phi}_r:\Omega\to M$ by
\[
\widetilde{\Phi}_r(x)=\widetilde{\exp}_x(r\widetilde{\nabla} u(x)).
\]
We have, for $x\in U$
\[
\widetilde{\Delta}_{{f_n}} u(x)\le n h(x)^{\frac{1}{n-1}};
\]
see \cite{Johne-arxiv} or \cite[Lemma 2.1]{IRV}.
Fix $\eta\in (0,1)$, and let $\tilde T$ be the radius given by Lemma \ref{lem:Lip}.
Define the sets
\[
\mathcal D_r(\Omega):=\{p \in M\,:\, \widetilde{\dist}(x, p) < r\text{ for all }x \in \Omega\}
\]
and
\[
\hat {\mathcal D}_{r}(\Omega):=\{p\in D_r(\Omega)\,:\,p=\sigma(t)\text{ for some }\sigma\in \mathcal{R}_0\text{ s.t. } [\sigma]\in V, t>\max\{\tilde{T},t_{*}\}\}.
\]

We are going to show that $\hat {\mathcal D}_{r}(\Omega)$
is contained in the image $\widetilde{\Phi}_r(A_r)$ of $A_r$ through the transport map.
\begin{lemma}\label{lem:image2}
    Let $p\in \hat{\mathcal D}_{r}(\Omega)$. Then, $p=\widetilde{\Phi}_r(y)$ for some $y\in A_r$. Moreover, the geodesic $s\mapsto \widetilde{\exp}_y(s\widetilde{\nabla}u)$, $s\in [0,r]$ does not intersect $B_{t_1}(o)$.
\end{lemma}
\begin{proof}
The existence of $y$ is a consequence of \cite[Lemma 2.2]{brendle}.
An inspection of the proof of that lemma shows that the point $y$ has the following characterization:
\begin{quote}    
the path $s\in[0,r]\mapsto \widetilde{\exp}_y(s\widetilde{\nabla} u(y))$
 minimizes the functional $\widetilde{\mathcal F}_u(\gamma):=u(\gamma(0))+\frac{1}{2}
\int_0^r |\gamma'(t)|_{\tilde g}^2\,dt$
among all smooth paths 
$\gamma:[0,r]\to M$ satisfying 
$\gamma(0)\in\Omega$ and $\gamma(r)=p$. As the minimizer is a constant speed minimizing geodesic, it holds in particular that the minimum of $\widetilde{\mathcal F}_u$ is $u(y)+\frac{1}{2r} \widetilde{\dist}(y,p)^2$.
\end{quote}
As $u(y)\ge u(x_\ast)$ by definition of $x_\ast$, necessarily $\widetilde{\dist}(y,p)\le \widetilde{\dist}(x_\ast,p)$. Recall that $p=\sigma(t)$ for some $\sigma\in \mathcal{R}_0$ such that $[\sigma]\in V$. Since the minimizing geodesic from $y$ to $p=\sigma(t)$ has length shorter than $\widetilde{\dist}(x_\ast,\sigma(t))$, Lemma \ref{lem:distances} allows to conclude that necessarily it can not intersect $B_{t_1}(o)$.
\end{proof}

Now, recall that $\Delta\varphi+\lambda (n-1)(\kappa\circ r_o) \varphi \le 0$ on $M\setminus B_{t_1}(o)$. By a straightforward computation the weighted manifold $(M\setminus B_{t_1}(o),\tilde g ,e^{-{f_n}}d\widetilde{\mathrm{vol}})$ satisfies
$\widetilde{\mathrm{Ric}}_{{f_n}}\ge 0$; see the proof of \cite[Proposition 2.3]{IRV} for details. Accordingly, for any $y\in A_r$ such that $p=\widetilde{\Phi}_r(y)\in\hat{\mathcal D}_{r}(\Omega)$, we can compute $|\det D\widetilde{\Phi}_r|(y)$ as in \cite[Proposition 3.3]{IRV}, obtaining that
\begin{align*}
|\det D\widetilde{\Phi}_r|(y)&\le e^{(2-n)({f}(p)-{f}(y))}\left[1 
+h(y)^{\frac{1}{n-1}}\int_0^r e^{\frac{2}{n}(2-n)\left({f}(\widetilde{\exp}_y(s\widetilde{\nabla}u(y))-{f}(y))\right)}\,ds
\right]^n\\
&\le e^{(2-n)({f}(p)-{f}(y))}\left[1+h(y)^{\frac{1}{n-1}}\int_0^r e^{\frac{2}{n}(2-n)\left({f}(\widetilde{\exp}_y(s\widetilde{\nabla}u(y)))\right)}\,ds
\right]^n.
\end{align*}
Recall that $(\log (1-t_1^{-(n-2+\epsilon)/2} - t_1^{-2}))/\lambda \le f(x) \le 0$ and set  \[
k:=\sup_{M\setminus B_{t_1}(o)}|f_n|=\frac{(n-2)}{\lambda}  \left|\log (1-t_1^{-(n-2+\epsilon)/2} - t_1^{-2})\right|<\infty.\]
Following again \cite{IRV}, we can write
\begin{align*}
|\mathrm{det} D\widetilde{\Phi}_r|(y)e^{-f_n(\widetilde{\Phi}_r(y))} 
&\le e^{2k}\left[ \frac{e^{-2k/n}}{r} +h(y)^{\frac{1}{n-1}}\right]^n r^n e^{-f_n(y)} .
\end{align*}
By Lemma \ref{lem:image2} and the area formula, 
\begin{align*}
 \int_{\hat{\mathcal D}_r(\Omega)}  e^{-f_n}d\widetilde{\mathrm{vol}}&\le \int_{\widetilde{\Phi}_r^{-1}(\hat{\mathcal D}_r(\Omega))} |\det D\widetilde{\Phi}_r| e^{-f_n\circ \widetilde{\Phi}_r}d\widetilde{\mathrm{vol}}
\\
&\le \int_\Omega e^{2k}\left[ \frac{e^{-2k/n}}{r} +h^{\frac{1}{n-1}}\right]^n r^n e^{-f_n}d\widetilde{\mathrm{vol}}.
\end{align*}

We are going to use the following result, whose proof is postponed.
\begin{lemma}\label{lemma:AVR}
It holds
\[
\limsup_{r\to\infty} \frac{\vol(\hat{\mathcal D}_{r}(\Omega))}{r^n}\ge \frac{(1-\eta)^{2(n-1)}}{n}\mu(\alpha,\theta,n).\]
\end{lemma}

Accordingly, for $n\ge 3$,
\begin{align*}
e^{2k}\int_\Omega h^\frac{n}{n-1}e^{-f_n}d\widetilde{\mathrm{vol}}&=e^{2k}\limsup_{r\to\infty}\int_\Omega \left[ \frac{e^{-2k/n}}{r} +h^{\frac{1}{n-1}}\right]^n e^{-f_n}d\widetilde{\mathrm{vol}}\\
&\geq \limsup_{r\to\infty}r^{-n}\int_{\hat{\mathcal D}_{r}(\Omega)}  e^{2f}d\mathrm{vol}\\
&\geq e^{-\frac{2k}{n-2}}\limsup_{r\to\infty}r^{-n}\vol(\hat{\mathcal D}_{r}(\Omega))\\
&\geq e^{-\frac{2k}{n-2}}\frac{(1-\eta)^{2(n-1)}}{n}\mu(\alpha,\theta,n),
\end{align*}
from which we deduce that
\begin{align*}
\left(\frac{(1-\eta)^{2(n-1)}}{n}\mu(\alpha,\theta,n)\right)^{\frac{1}{n}}&e^{\frac{-(n-1)k}{n}}\left(\int_\Omega \ h^{\frac{n}{n-1}}d\mathrm{vol}\right)^{\frac{n-1}{n}}\\&\leq \left(\frac{(1-\eta)^{2(n-1)}}{n}\mu(\alpha,\theta,n)\right)^{\frac{1}{n}}\left(\int_\Omega \ h^{\frac{n}{n-1}}e^{-f_n}d\widetilde{\mathrm{vol}}\right)^{\frac{n-1}{n}}\\
&\le e^{\frac{2(n-1)k}{n(n-2)}} \int_\Omega \ h^{\frac{n}{n-1}}e^{-f_n}d\widetilde{\mathrm{vol}} \\
&=\frac{e^{\frac{2(n-1)k}{n(n-2)}}}{n}\left(\int_\Omega |\widetilde{\nabla h}|_{\tilde g}e^{-f_n}\,d\widetilde{\mathrm{vol}} + \int_{\partial \Omega} h e^{-f_n}\,d\widetilde{\mathrm{vol}}_{n-1}\right)\\
&= \frac{e^{\frac{2(n-1)k}{n(n-2)}}}{n}\left(\int_\Omega |\nabla h|\,d\mathrm{vol} + \int_{\partial \Omega} h e^{f}\,d\mathrm{vol}_{n-1}\right)\\
&\le \frac{e^{\frac{2(n-1)k}{n(n-2)}}}{n}\left(\int_\Omega |\nabla h|\,d\mathrm{vol} + \int_{\partial \Omega} h \,d\mathrm{vol}_{n-1}\right),
\end{align*}
the first equality being due to the scaling assumption \eqref{e:scalingconf}. Summarizing, we have proved that on any $\Omega\subset M\backslash B_{T_0}(o)$ it holds
\[
\left(\int_\Omega \ h^{\frac{n}{n-1}}d\mathrm{vol}\right)^{\frac{n-1}{n}}\leq \frac{\left(\frac{(1-\eta)^{2(n-1)}}{n}\mu(\alpha,\theta,n)\right)^{-\frac{1}{n}}e^{\frac{(n-1)k}{n-2}}}{n}\left(\int_\Omega |\nabla h|\,d\mathrm{vol} + \int_{\partial \Omega} h \,d\mathrm{vol}_{n-1}\right).
\]
Similarly, when $n=2$ we obtain the constant $\frac{1}{2}(\frac{(1-\eta)^2}{2}\mu(\alpha,\theta,2))^{-\frac{1}{2}}e^{2\sup|f|}$. 

Note that in the final estimate $\eta$ can be arbritrarily chosen in $(0,1)$ and the other constants are independent of $\eta$. In particular, one can let $\eta$ go to $1$ and get rid of the dependence on $\eta$ in the previous inequality.
\end{proof}

It remains only to prove Lemma \ref{lemma:AVR}.

\begin{proof}[Proof (of Lemma \ref{lemma:AVR})]
According to Lemma \ref{lem:Lip}, for all $\sigma,\gamma\in\mathcal R_o$,
\begin{equation}\label{e:LipLoc}
    \frac{1}{t}\dist(\sigma(t),\gamma(t))\ge  (1-\eta)^2\angle_\infty([\sigma],[\gamma]).
\end{equation}
as soon as $t>\tilde T$ and $\dist(\sigma(t),\gamma(t))<C_\eta\,\tilde T$. Here, $\eta$ and $\tilde T$ are the constant defined before Lemma \ref{lem:image2}.

Now, note that $\widetilde{\dist}\le \dist$ so that $\mathcal D_r(\Omega)\supset \widetilde{B_{r-\tilde r_\Omega}}(o)\supset B_{r-\tilde r_\Omega}(o)$, where $\tilde r_\Omega$ is the smaller radius such that $\widetilde {B_{\tilde r_\Omega}}(o)\supset\Omega$. Accordingly.
\[
\hat {\mathcal D}_{r}(\Omega)\supset\{\sigma(t)\,:\, \max\{\tilde T;t_{*}\} < t < r-r_\Omega\text{ and }\sigma\in \mathcal{R}_0\text{ s.t. } [\sigma]\in V\},
\]
so that, letting $T_{\Omega}:=\max\{\tilde{T}; t_{*}\}$,
\begin{align*}
        \mathcal H^{n}(\hat  {\mathcal D}_{r}(\Omega))&\ge \int_{T_{\Omega}}^{r-r_\Omega} \mathcal H^{n-1}(\hat  {\mathcal D}_{r}(\Omega)\cap S_t)\,dt \\
    &\ge (1-\eta)^{2(n-1)}\mathcal H^{n-1}(V)\int_{T_{\Omega}}^{r-r_\Omega} t^{n-1}\,dt,
\end{align*}
where we used the coarea formula in the first inequality \cite{Fe}, and the fact that the map  from $\hat  {\mathcal D}_{r}(\Omega)\cap S_t$ to $(V,\angle_{Tits})$ defined by $\sigma(t)\mapsto [\sigma]$ is locally $\frac{1}{(1-\eta)^{2}t}$-Lipschitz due to \eqref{e:LipLoc}; see for instance \cite[Proposition 1.7.8]{BBI}.

To conclude, from \eqref{e:vol} we note that 
\[
\mathcal H^{n-1}(V)\ge \mu(\alpha,\theta,n).\]
All together, we have proved that
  \[
  \mathcal H^{n}(\hat  {\mathcal D}_{r}(\Omega))\ge
\frac{(1-\eta)^{2(n-1)}}{n} \mu(\alpha,\theta,n)[(r-r_\Omega)^n-T_{\Omega}^n].
\]
  \end{proof}

\begin{remark}\label{Rmk:ToyCases}{\rm
In the stronger assumptions of Corollary \ref{coro:nonnegoutcomp}, that is when the sectional curvature of $M$ is non-negative outside some compact set $K$ of $M$, there is no need to conformally deform the metric of $M$. In particular, one can replace $\varphi\equiv 1$ and $f\equiv 0$ in the above proof.

An even simpler case is when  $M$ has nonnegative sectional curvatures outside a compact set $K\subset M$, $M$ has only one end, and $\diam M_\infty<\pi$. Indeed in this case, 
we \textsl{claim} that there exists $R_0$ large enough so that any minimizing geodesic connecting two points of  $M\setminus B_{R_0}(o)$ does not intersect $K$. One can thus implement Brendle's approach to prove the isoperimetric inequality on $M\setminus B_{R_0}(o)$. Namely, it is enough to consider $\hat A_{r}:=\{y\in A_r\ :\ \Phi_r(y)\not\in B_{R_0}\}$, so that the minimizing geodesics from $y\in \hat A_{r}$ to $\Phi_r(y)$ is totally contained in $M\setminus K$, and the curvature is non-negative along that geodesic. Since removing $B_R(o)$ from the image of $\Phi_r$ does not affect the asymptotic volume growth, one can conclude the proof.

To show the above claim, suppose there exist two sequences of points $\{x_j\}$, $\{y_j\}$ such that $\dist(o,x_j)>j,\dist(o,y_j)>j$ and a unit-speed minimizing geodesic $\gamma_j$ from $x_j$ to $y_j$ which intersects $K$.
On the one hand, $\dist(x_j,y_j)\ge 2j-2\diam{K}$. On the other hand $\lim_{j\to\infty} \dist(x_j,y_j)/j<2$ by Lemma \ref{lem:mno} (2) and the assumption on $M_\infty$, thus giving a contradiction. %Note that in particular this implies that there are no lines in $M$ when $\diam M_\infty<\pi$. 

A similar argument permits also to deal the case in which there are several ends, but each connected component of $M_\infty$ has diameter $<\pi$. 
}
\end{remark}

\begin{remark}\label{rmk:constants}
    {\rm
According to the previous remark, when all the connected components of $M_\infty$ have diameter $<\pi$ the isoperimetric constant in \eqref{e:isoper asymp} restricted to $C^\infty_c(E\setminus B_{R_0}(o))$, where $E$ is an end of $M$, is the sharp Euclidean one, i.e. $n\omega_n^{1/n}\theta^{1/n}$ with $\theta=\limsup_{r\to\infty} \omega_n^{-1}r^{-n}\vol(B_r(0)\cap E)$ the AVR of $E$. Clearly, the constant in \eqref{e:isoper asymp} for general functions in $C^\infty_c(M)$ necessarily depends on the whole manifold and can not be estimated only in term of the quantities appearing in our assumptions. For instance, adding a tubular handle $[0,L]\times \epsilon\mathbb S^{n-1}$ to $M$ near $o$ makes the isoperimetric constant explode when $L\to\infty$ without adding negative curvature outside a compact neighborhood of $o$. 
    }
\end{remark}
   
\section{Submanifolds in ambient manifolds with nonnegative sectional curvature outside a compact set}\label{sec:Michael-Simon}

Throughout this section we will assume that $(M^{n+m},g)$ is a complete smooth $(n+m)$-dimensional Riemannian manifold satisfying
\[
\mathrm{Sect}^M\ge 0
\]
outside some compact set $K\subset M$. Let $N=n+m$. We will further assume that, given a reference point $o\in M$, $\textrm{AVR}>0$ on each end $E_j$ of $M$, i.e. there exists $\theta>0$ such that
    \[
    \limsup _{r\to\infty}\frac{\vol(B^M_r(o)\cap E_j)}{r^{N}}\geq \theta>0, \qquad \forall\,E_{j}.
    \] 
Let $\Sigma$ be a compact submanifold of dimension $n$. We will denote by $\overline{\nabla}$ the Levi-Civita connection on $M$ and by $\nabla^\Sigma$ the induced connection on $\Sigma$. The second fundamental form $II$ of $\Sigma$ is given by
\[
g( II(X,Y),N)=-g( \overline{\nabla}_X N, Y),
\]
where $X,\, Y$ are tangent vector fields on $\Sigma$ and $N$ is a normal vector field along $\Sigma$. Moreover, the mean curvature vector of $\Sigma$ is defined by $H=\mathrm{tr}II$.

Since every end has positive AVR, Section \ref{sec:preliminar} gives a positive constant $\mu_{\alpha,N}:=\mu(\alpha, \theta, N)>0$ such that, for every ray $\sigma\in\mathcal{R}_{o}$,
\[
\ \mathcal{H}^{N-1}(B_{\alpha}[\sigma])\geq\mu_{\alpha,N}.
\]

For the reader's convenience, we recall here some objects and notations that were already introduced in Section \ref{sec:asymptotic} and will be used throughout the whole section. We fix $\alpha\in(0,\frac{\pi}{3})$ as in Section \ref{sec:preliminar} and choose $\delta>0$ so small that
\begin{equation}\label{eq: 5.1}
\ 2\sin\frac{\alpha}{2}<1-2\delta
\end{equation}
Choose $a, \epsilon\in(0,1)$ such that $1+\epsilon-(1-\epsilon)\sqrt{1-a^2}<\delta/2$. Let $\alpha_{0}>0$ be the constant given by Lemma \ref{lem:topo2} and let $\rho_{0}$ be the radius in Lemma \ref{lem:density radial} applied with $\zeta=\alpha_{0}$. Let $T=T(\delta)$ be the radius in Lemma \ref{lem:unif lim}, so that
\[
\ \frac{\mathrm{dist}(\sigma_{1}(t),\sigma_{2}(t))}{t}\leq2\sin\left(\frac{\angle_{\infty}([\sigma_{1}], [\sigma_{2}])}{2}\right)+\delta
\]
for every $t\geq T$ and every pair of rays $\sigma_{1},\sigma_{2}\in\mathcal{R}_{o}$. We shall also fix a $\eta$ in Lemma \ref{lem:Lip} and let $\tilde{T}$ be the corresponding radius in Lemma \ref{lem:Lip}, and choose $t_{0}\gg1$ so that $K\subset B_{t_{0}/2}(o)$ and $t_{0}\geq\max\left\{\rho_{0},T,\tilde{T}\right\}$.

\subsection{Michael-Simon inequality}
We are going to prove that the Michael-Simon-type inequality \eqref{eq:MichaelSimon} holds on any compact submanifold $\Sigma$ in $M\setminus B_{t_0}(o)$. Note that it suffices to prove the result in case $\Sigma$ is connected (we refer e.g. to \cite{brendle-Rn} for the discussion of the general case). Let $h$ be a positive smooth function on $\Sigma$. We assume by scaling that
\begin{equation}\label{e:scalingsub}   
\int_\Sigma \sqrt{|\nabla^\Sigma h|^2+h^2|H|^2}\,d\mathrm{vol}_\Sigma + \int_{\partial \Sigma} h \,d\mathrm{vol}_{n-1} = n \int_\Sigma h^{\frac{n}{n-1}}\,d\mathrm{vol}_\Sigma,
\end{equation}
and we consider the Neumann problem
\[
\begin{cases}
\mathrm{div}^\Sigma(h\nabla^\Sigma u)= nh^{\frac{n}{n-1}}-\sqrt{|\nabla^\Sigma h|^2+h^2|H|^2},&\textrm{in }\Sigma\\
\langle \nu,\nabla^\Sigma u\rangle = 1,& \text{on }\partial\Sigma.
\end{cases}
\]
The solution is unique up to an additive constant and continuous on $\bar\Sigma$. As in \cite{brendle}, define
 \[
 U=\{x\in\Sigma\setminus\partial\Sigma\, ,y\in T^\perp_x\Sigma:\,|\nabla^\Sigma u(x)|^2+|y|^2<1\}
 \]
 and, for $r>0$,
 \begin{equation*}
A_r=\{(\bar{x},\bar{y})\in U\,:\,\forall (x,y)\in U,\ ru(x)+\frac{1}{2}d(x,\exp_{\bar{x}}(r\nabla^\Sigma u(\bar{x})+r\bar{y}))^2 \ge ru(\bar{x}) +\frac{1}{2}r^2(|\nabla^\Sigma u(\bar{x})|^2+|\bar{y}|^2)\} 
 \end{equation*}
Define the $C^{1,\alpha}$
transport map $\Phi_r:T^\perp\Sigma\to M$ by
\[
\Phi_r(x,y)=\exp_x(r\nabla^\Sigma u(x)+ry),
\]
with $\exp$ the exponential map of $M$. The following lemma is due to Brendle \cite{brendle} and its proof is independent of the curvature condition of the ambient space. 

\begin{lemma}[Lemma 4.2 in \cite{brendle}]\label{Lemma4.2}
For each $0\leq\mu<1$, the set
\[
\mathcal B_{r,\mu}(\Sigma):=\{p \in M\,:\, \mu r< \dist(x, p) < r\text{ for all }x \in \Sigma\}
\]
is contained in $\Phi_r(\{(x,y)\in A_r\,:\, |\nabla^\Sigma u(x)|^2+|y|^2>\mu^2\})$.
\end{lemma}

Let $x_{*}\in\Sigma$ be a point where $u$ attains its minimum and write $x_{*}=\sigma_{*}(t_{*})$, where $\sigma_{*}$ is a minimizing geodesic from $o$ to $x_{*}$. Since $\Sigma\subset M\setminus B_{t_{0}}(o)$, we have $t_{*}\geq t_{0}\geq \rho_{0}$. Hence, by Lemma \ref{lem:density radial}, there exists a ray $\gamma_{*}\in\mathcal{R}_{o}$ such that $\angle_{o}(\sigma_{*},\gamma_{*})<\alpha_{0}$. Set now $V\doteq B_{\alpha}([\gamma_{*}])\subset M_{\infty}$ and let $\hat{\mathcal{B}}_{r,\mu}(\Sigma)$ be the subset of $\mathcal{B}_{r,\mu}(\Sigma)$ defined by
\[
\ \hat{B}_{r,\mu}(\Sigma):=\left\{p\in\mathcal{B}_{r,\mu}(\Sigma)\,:\, p=\sigma(t)\, \text{for some }\sigma\in\mathcal{R}_{o},\,[\sigma]\in V, t>t_{*}\right\}.
\]
\begin{lemma}\label{lem:image_sub} Let $p=\sigma(t)\in \hat {\mathcal B}_{r,\mu}(\Sigma)$ and let $(x,y)\in A_r$ be such that $\Phi_r(x,y)=p$. Then the geodesic $s\in[0,r]\mapsto \exp_x(s\nabla^\Sigma u(x)+sy)$ is contained in $M\setminus B_{t_0/2}(o)\subset  M\setminus K$.
\end{lemma}
\begin{proof}
First we prove that for any  $p=\Phi_r(x,y)\in \mathcal B_{r,\mu}(\Sigma)$ with $(x,y)\in A_r$, it holds $\dist(x,p)\le \dist(x_\ast,p)$, where $x_\ast\in\Sigma$ is a point where $u$ attains its minimum. Indeed, an inspection of the proof of \cite[Lemma 4.2]{brendle} shows that $(x,y)$ has the following characterization:
\begin{quote}    
the path $s\in[0,r]\mapsto \exp_x(s\nabla^\Sigma u(x)+sy)$
 minimizes the functional $\mathcal F_u(\gamma):=u(\gamma(0))+\frac{1}{2}
\int_0^r |\gamma'(t)|^2\,dt$
among all smooth paths 
$\gamma:[0,r]\to M$ satisfying 
$\gamma(0)\in\Sigma$ and $\gamma(r)=p$. As the minimizer is a constant speed geodesic, it holds in particular that the minimum of $\mathcal F_u$ is $u(x)+\frac{1}{2r} \dist(x,p)^2$.
\end{quote}
As $u(x)\ge u(x_\ast)$ by definition of $x_\ast$, necessarily $\dist(x,p)\le \dist(x_\ast,p)$.

On the other hand, reasoning as in the first part of the proof of Lemma \ref{lem:distances} with the choice $f\equiv 0$, it can be proved that for all $p=\sigma(t)\in \hat {\mathcal B}_{r,\mu}(\Sigma)$,  it holds $\dist(p,x_\ast)< t=r_o(p)$. Accordingly, if $(x,y)\in A_r$ is such that $\Phi_r(x,y)=p$, then $\dist(x,p)<r_o(p)$. In particular, by the triangle inequality, the minimizing geodesic $s\in[0,r]\mapsto \check\gamma(s):=\exp_x(s\nabla^\Sigma u(x)+sy)$ is contained in $M\setminus B_{t_0/2}(o)\subset  M\setminus K$. Indeed, suppose by contradiction that $\check\gamma(s)\in B_{t_0/2}(o)$ for some $s\in[0,r]$. Then
\begin{align*}
    r_0(p)\ge \dist(x,p)&= \dist(x,\check\gamma(s))+\dist(p,\check\gamma(s))\\
    &> \dist(x,B_{t_0/2}(o))+ \dist(p,B_{t_0/2}(o))\\
    &\ge (t_0-t_0/2)+(r_o(p)-t_0/2) = r_0(p),
\end{align*}
giving a contradiction.
\end{proof}

A key point in Brendle's proof of the Michael-Simon-type inequality for submanifolds of a manifold with nonnegative sectional curvature   is the estimate of the Jacobian of $\Phi_r$; see \cite[Proposition 4.8]{brendle}. An inspection of the proof shows that the aforementioned estimate holds whenever $M$ has non negative sectional curvature outside $B^M_{R_0}(o)$ for some $R_0$ and
$x$, $\Phi_r(x,y)$ and $\exp_x(s\nabla^\Sigma u(x)+sy)$ all belong to $M\backslash B^M_{R_0}(o)$
for any $(x,y)\in A_r$, $s\in[0,r]$. 

We thus have the validity of the following

\begin{proposition}\label{Prop4.8general}
For all $(x,y)\in A_r$ such that $\Phi_r(x,y)\in \hat{\mathcal B}_{r,\mu}(\Sigma)$ it holds    
\[
|\mathrm{det} D\Phi_r|( x, y) \le r^m\left(1+h( x)^{\frac{1}{n-1}}r\right)^n.
\]   
\end{proposition}

With this preparation we can now complete the proof of Theorem \ref{thm:nonnegsuboutcomp}.
\begin{proof}[Proof of Theorem \ref{thm:nonnegsuboutcomp}]
We first focus on the case $m\geq 2$. Using the coarea formula together with Proposition \ref{Prop4.8general} and computing as in \cite[p.20]{brendle}, we deduce that 
\begin{align*}
\mathrm{vol}(\hat{\mathcal B}_{r,\mu}(\Sigma))&\leq \int_{\{x\in\Sigma: |\nabla^\Sigma u(x)|<1\}} \left(\int_{\{y\in T_x^\perp\Sigma: \mu^2<|\nabla^\Sigma u(x)|^2+|y|^2<1\}}|\mathrm{det}D\Phi_r(x,y)|dy\right)d\mathrm{vol}_\Sigma(x)\\
&\leq \frac{m}{2}\vol(\mathbb{B}^m_1(0))(1-\mu^2)\int_\Sigma  r^{m+n}\left(\frac{1}{r}+h(x)^{\frac{1}{n-1}}\right)^n d\mathrm{vol}_\Sigma(x).
\end{align*}
Dividing by $r^{n+m}$ and sending $r\to\infty$ while keeping $\mu$ fixed we get 
\[
\limsup_{r\to\infty} \frac{\vol(\hat{\mathcal B}_{r,\mu}(\Sigma))}{r^{n+m}}\leq \frac m2\vol(\mathbb{B}^m_1(0))(1-\mu^2)\int_\Sigma h^{\frac{n}{n-1}},
\]
for all $0\leq\mu<1$. Now, let $r_\Sigma$ be the smaller radius such that $ {B_{r_\Sigma}}(o)\supset\Sigma$. For $r$ large enough, 
\[
\hat {\mathcal B}_{r,\mu}(\Sigma)\supset\{\sigma(t)\,:\, \max\{t_{\ast};\mu r+r_\Sigma\} < t < r-r_\Sigma\text{ and }\sigma\in \mathcal{R}_0\text{ s.t. } [\sigma]\in V\}.
\]
By Lemma \ref{lem:Lip}, the projection $\pi_{t}:\left\{\sigma(t):\,[\sigma]\in V\right\}\to V$, $\pi_{t}(\sigma(t))=[\sigma]$ is locally $((1-\eta)^2t)^{-1}$-Lipschitz for $t>\tilde{T}$. Hence
\[
\ \mathcal{H}^{N-1}\left(\left\{\sigma(t)\,:\, [\sigma]\in V\right\}\right)\geq (1-\eta)^{2(N-1)}t^{N-1}\mathcal{H}^{N-1}(V)\geq (1-\eta)^{2(N-1)}\mu_{\alpha, N}t^{N-1}.
\]
It follows that
\[
\ \mathrm{vol}(\hat{\mathcal{B}}_{r,\mu}(\Sigma))\geq (1-\eta)^{2(N-1)}\mu_{\alpha, N}\int_{\max\left\{t_{*};\mu r+r_{\Sigma}\right\}}^{r-r_{\Sigma}}t^{N-1}dt.
\]
Consequently,
\begin{align*}
    \limsup_{r\to\infty} \frac{\vol(\hat{\mathcal B}_{r,\mu}(\Sigma))}{r^{N}}&\ge \frac{\Theta_{MS}}{N} (1-\mu^{N}),
\end{align*}
where $\Theta_{MS}\doteq(1-\eta)^{2(N-1)}\mu_{\alpha,N}$.
Hence,
\[
\frac{\Theta_{MS}}{N} (1-\mu^{N})\leq \frac m2\vol(\mathbb{B}^m_1(0))(1-\mu^2)\int_\Sigma h^{\frac{n}{n-1}}.
\]
Dividing both sides by $1-\mu$ and letting $\mu\to1$, we thus obtain 
\[
\Theta_{MS}\leq m\vol(\mathbb{B}^m_1(0))\int_\Sigma h^{\frac{n}{n-1}}.
\]
Finally, by the initial rescaling assumption,
\begin{align*}
\int_\Sigma \sqrt{|\nabla^\Sigma h|^2+h^2|H|^2}\,d\mathrm{vol}_\Sigma + \int_{\partial \Sigma} h \,d\mathrm{vol}_{n-1} &= n \int_\Sigma h^{\frac{n}{n-1}}\,d\mathrm{vol}_\Sigma\\
&=n \left(\int_\Sigma h^{\frac{n}{n-1}}\,d\mathrm{vol}_\Sigma\right)^{\frac{1}{n}}\left(\int_\Sigma h^{\frac{n}{n-1}}\,d\mathrm{vol}_\Sigma\right)^\frac{n-1}{n}\\
&\geq n\left(\frac{1}{m \vol(\mathbb{B}^m_1(0))}\right)^{\frac{1}{n}}\Theta_{MS}^\frac{1}{n} \left(\int_\Sigma h^{\frac{n}{n-1}}\,d\mathrm{vol}_\Sigma\right)^{\frac{n-1}{n}}
\end{align*}
and this conclude the proof of the validity of equation \eqref{eq:MichaelSimon} in Theorem \ref{thm:nonnegsuboutcomp}. 

It remains to consider the case $m=1$. In this case the annular estimate used above is no longer effective. We use instead the weighted area formula argument appearing in \cite{Brendle_Survey}. Set $N=n+1$, $\mu=0$ and
\[
\ \mathcal{I}_{r}:=\int_{A_{r}}\frac{|\mathrm{det}D\Phi_{r}(x,y)|}{\sqrt{1-|\nabla^{\Sigma}u(x)|^2-|y|^2}}dyd\mathrm{vol}_{\Sigma}(x).
\]
Then we get that
\[
\ \mathcal{I}_{r}\leq \pi r^{N}\int_{\Sigma}\left(\frac{1}{r}+h^{\frac{1}{n-1}}\right)d\mathrm{vol}_{\Sigma},
\]
and, with the same notations as above,
\begin{equation*}
\mathcal{I}_{r}\geq \int_{\hat{\mathcal{B}}_{r,0}(\Sigma)}\frac{d\mathrm{vol}_{M}(p)}{\sqrt{1-\left(\frac{r_{o}(p)-r_{\Sigma}}{r}\right)^2}}\geq \int_{t^{*}}^{r-r_{\Sigma}}\frac{\mathcal{H}^{N-1}\left(\left\{\sigma(t)\,:\,\left[\sigma\right]\in V\right\}\right)}{\sqrt{1-\left(\frac{t-r_{\Sigma}}{r}\right)^2}}dt.
\end{equation*}
By Lemma \ref{lem:Lip}, we hence get that
\[
\mathcal{I}_{r}\geq \Theta_{MS}\int_{t^{*}}^{r-r_{\Sigma}}\frac{t^{N-1}}{\sqrt{1-\left(\frac{t-r_{\Sigma}}{r}\right)^2}}dt.
\]
Reasoning now as in the previous case, we eventually get that for $m=1$ the constant in the Michael-Simon inequality is 
\[
\ n\left(\frac{\Theta_{MS}\int_{0}^{1}\frac{s^{N}}{\sqrt{1-s^{2}}}ds}{\pi}\right)^{\frac{1}{n}}.
\]

As for equation \eqref{eq:MichaelSimonunbounded}, note that if $\varphi:\Sigma^n\mapsto M^{n+m}$ is a proper immersion of $\Sigma$ in $M$, the set $D=\varphi^{-1}(B_{R_0}(o))$ is compact and hence the Michael-Simon-type inequality holds on $\Sigma\backslash D$. We can now apply the generalization of Theorem \ref{th:extension} pointed out in Remark \ref{rmk:extension H}
to get the validity of the Michael-Simon-type inequality on the whole $\Sigma$.
\end{proof}

\begin{remark}{\rm
Note that if we further assume that the ambient manifold $M$ is such that every connected component of $M_\infty$ has diameter smaller than $\pi$, then the sharp Brendle's Michael-Simons inequality holds outside a suitable compact set. Indeed in this case, the reasoning in Remark \ref{Rmk:ToyCases} implies that there exists $R_0$ large enough so that any minimizing geodesic connecting two points of the same connected component of $M\setminus B_{R_0}(o)$ does not intersect $K$. One can thus implement Brendle's approach to prove the isoperimetric inequality on $M\setminus B_{R_0}(o)$. To this end, as in the case of subdomains, it is enough to consider $\hat A_{r}:=\{(x,y)\in A_r\ :\ \Phi_r(x,y)\not\in B_{R_0}\}$, so that the minimizing geodesics from $x\in \hat A_{r}$ to $\Phi_r(x,y)$ is totally contained in $M\setminus K$, and the curvature is non-negative along that geodesic. Since removing $B_R(o)$ from the image of $\Phi_r$ does not affect the asymptotic volume growth, one can conclude the proof.   
}
\end{remark}

\subsection{Log-Sobolev inequality} 
The technique used in the previous subsection can also be used to prove the validity of a Log-Sobolev inequality for compact submanifolds, thus generalizing previous analogous results obtained by S. Brendle in the case of submanifolds of the Euclidean space \cite{brendle-LogSob}, and by  C. Yi and Y. Zheng when the curvature of the ambient space is everywhere non-negative \cite{YZ}. More precisely   \begin{theorem}\label{thm:LogSob}
    Suppose that $(M^{n+m},g)$ has nonnegative sectional curvatures outside a compact set $K\subset M$. Let $o\in K$. Suppose that 
    $\textrm{AVR}>0$ on each end $E_j$ of $M$, i.e. there exists $\theta>0$ such that
    \[
    \limsup _{r\to\infty}\frac{\vol(B^M_r(o)\cap E_j)}{r^{n+m}}\geq \theta>0,\qquad\forall\,E_{j}.
    \] 
    Then, there exists $t_0>0$ and a positive constant $\Theta_{LS}$, depending only on $n$, $m$ and the uniform lower bound $\theta$, such that for any compact $n$-dimensional submanifold without boundary $\Sigma\subset M\backslash B_{t_0}(o)$ and any positive smooth function $h$ on $\Sigma$ it holds
    \begin{align}\label{eq:LogSob}
    &\int_\Sigma h\left(\log h+ n+\log\Theta_{LS}+\frac{n}{2}\log(4\pi)\right) d\mathrm{vol}_\Sigma -\int_\Sigma \frac{|\nabla^\Sigma h|^2}{h}\,d\mathrm{vol}_\Sigma-\int_\Sigma h|H|^2\,d\mathrm{vol}_\Sigma\\\nonumber
&\leq \left(\int_\Sigma h d\mathrm{vol}_\Sigma\right)\log\left(\int_\Sigma  h d\mathrm{vol}_\Sigma\right).
    \end{align}
\end{theorem}
We are going now to prove that the Log-Sobolev inequality \eqref{eq:LogSob} holds on any compact submanifold $\Sigma$ without boundary in $M\setminus B_{t_0}(o)$. Here we are using the constants fixed at the beginning of this section. Note that it suffices to prove the Theorem in case $\Sigma$ is connected (we refer e.g. to \cite{YZ} for the discussion of the general case). Let $h$ be a positive smooth function on $\Sigma$. By scaling, we may assume that
\begin{equation}\label{eq:normLogSob}
\int_\Sigma h\log h \,d\mathrm{vol}_\Sigma -\int_\Sigma \frac{|\nabla^\Sigma h|^2}{h}\,d\mathrm{vol}_\Sigma-\int_\Sigma h|H|^2\,d\mathrm{vol}_\Sigma=0.
\end{equation}
As a consequence of the previous assumption and of the fact that $\Sigma$ is connected, there exists a smooth function $u:\Sigma\to\mathbb{R}$ satisfying
\[
\mathrm{div}^\Sigma (h\nabla^\Sigma u)=h\log h-\frac{|\nabla^\Sigma h|^2}{h}-h|H|^2.
\]
For each $r>0$, let us define the set
 \begin{equation*}
C_r=\{(\bar{x},\bar{y})\in T^\perp\Sigma\,:\,\forall x\in \Sigma,\ ru(x)+\frac{1}{2}d(x,\exp_{\bar{x}}(r\nabla^\Sigma u(\bar{x})+r\bar{y}))^2 \ge ru(\bar{x}) +\frac{1}{2}r^2(|\nabla^\Sigma u(\bar{x})|^2+|\bar{y}|^2)\} 
 \end{equation*}

Let $\Phi_r$ be defined as in the previous subsection. The following lemma is due to Yi and Zheng \cite{YZ} (see also \cite{DDLb}) and its proof is independent of the curvature condition of the ambient space. 

\begin{lemma}[Lemma 3.2 in \cite{YZ}]\label{Lemma3.2}
For each $r>0$, $\Phi_r(C_r)=M$.
\end{lemma}

Let $x_{*}$ be a point where $u$ attains its minimum. Let $x_{*}=\sigma_{*}(t_{*})$ where $\sigma_{*}$ is a minimizing geodesic from $o$ to $x_{*}$. Since $\Sigma\subset M\setminus B_{t_{0}}(o)$ and we have $t_{*}\geq t_{0}\geq \max\{\rho_{0},T, \tilde{T}\}$, by Lemma \ref{lem:density radial} there exists $\gamma_{*}\in\mathcal{R}_{o}$ such that $\angle_{o}(\sigma_{*},\gamma_{*})<\alpha_{0}$. Set $V:=B_{\alpha}([\gamma_{*}])\subset M_{\infty}$.
Now consider the set
\[
\mathcal{D}=\{p\in M\,:\,p=\sigma(t),\,\sigma\in \mathcal{R}_0, [\sigma]\in V,\,t>t_{*}\}
\]
Reasoning as in the previous subsection one can show that, if $p=\sigma(t)\in \mathcal{D}$ and $(x,y)\in C_r$ are such that $\Phi_r(x,y)=p$, then the geodesic $s\in[0,r]\mapsto \exp_x(s\nabla^\Sigma u(x)+sy)$ is contained in $M\setminus B_{t_0/2}(o)\subset  M\setminus K$. We can then reason as in \cite{YZ} and prove the validity of the following 
\begin{proposition}\label{Prop4.8generalLogSob}
For all $(x,y)\in C_r$ such that $\Phi_r(x,y)\in \mathcal{D}$ it holds    
\[
e^{-\frac{d(x,\Phi_r(x,y))^2}{4r^2}}|\mathrm{det} D\Phi_r|( x, y) \le r^{n+m}h(x)e^{\frac{n}{r}-n}e^{-\frac{|2H(x)+y|^2}{4}}.
\]       
\end{proposition}

With this preparation we can now complete the proof of Theorem \ref{thm:LogSob}.
\begin{proof}[Proof of Theorem \ref{thm:LogSob}]
Using the coarea formula together with the triangle inequality and Proposition \ref{Prop4.8generalLogSob}, we deduce that 
\begin{align*}
\int_{\mathcal{D}} e^{-\frac{(d(p,o)+r_\Sigma)^2}{4r^2}}\,d\mathrm{vol}_M(p)&\leq \int_\Sigma \int_{ T_x^\perp\Sigma}e^{-\frac{d(x,\Phi_r(x,y))^2}{4r^2}}|\mathrm{det}D\Phi_r(x,y)|dy d\mathrm{vol}_\Sigma(x)\\
&\leq r^{n+m}e^{\frac{n}{r}-n}(4\pi)^{\frac{m}{2}}\int_\Sigma  h(x) d\mathrm{vol}_\Sigma(x),
\end{align*}
where $r_\Sigma$ is the smaller radius such that $ {B_{r_\Sigma}}(o)\supset\Sigma$.
Using Lemma \ref{lem:Lip} as in the previous subsection, for $t>t_{*}$ we have
\[
\ \mathcal{H}^{N-1}(\left\{\sigma(t)\,:\, [\sigma]\in V\right\})\geq (1-\eta)^{2(N-1)}\mu_{\alpha,N}t^{N-1}=\Theta_{MS}t^{N-1}.
\]
Therefore
\[
\ \int_{D}e^{-(d(o,p)+r_{\Sigma}^{2})/4r^2}d\mathrm{vol}_{M}(p)\geq\Theta_{MS}\int_{t_{*}}^{+\infty}t^{N-1}e^{-(t+r_{\Sigma}^2)/4r^2}dt.
\]
Dividing by $r^{N}$ and letting $r\to+\infty$, we get 
\[
\ \limsup_{r\to\infty}\frac{1}{r^{N}}\int_{D}e^{-(d(o,p)+r_{\Sigma})^{2}/4r^{2}}d\mathrm{vol_{M}(p)}\geq \Theta_{MS}\int_{0}^{+\infty}\tau^{N-1}e^{-\tau^2/4}d\tau.
\]
Letting 
\[
\ \Theta_{LS}:=\frac{\Theta_{MS}}{(4\pi)^{N/2}}\int_{0}^{\infty}\tau^{N-1}e^{-\tau^{2}/4}d\tau=\frac{(1-\eta)^{2(N-1)}\mu_{\alpha, N}}{(4\pi)^{N/2}}\int_{0}^{+\infty}\tau^{N-1}e^{-\tau^2/4}d\tau
\]
we hence get
\[
n+\log\Theta_{LS}+\frac{n}{2}\log(4\pi)\leq \log\left(\int_\Sigma  h d\mathrm{vol}_\Sigma\right).
\]
Combining this inequality with the normalization \eqref{eq:normLogSob}, we obtain
\begin{align*}
&\int_\Sigma h\left(\log h+ n+\log\Theta_{LS}+\frac{n}{2}\log(4\pi)\right) d\mathrm{vol}_\Sigma -\int_\Sigma \frac{|\nabla^\Sigma h|^2}{h}\,d\mathrm{vol}_\Sigma-\int_\Sigma h|H|^2\,d\mathrm{vol}_\Sigma\\
&=\int_\Sigma h\left(n+\log\Theta_{LS}+\frac{n}{2}\log(4\pi)\right) d\mathrm{vol}_\Sigma\\
&\leq \left(\int_\Sigma h d\mathrm{vol}_\Sigma\right)\log\left(\int_\Sigma  h d\mathrm{vol}_\Sigma\right).
\end{align*}
This completes the proof of Theorem \ref{thm:LogSob}.
\end{proof}
\bibliographystyle{alpha}
\bibliography{sample}

\end{document}